%% file: ocha62806cmp.tex
\documentclass[11pt]{article}
\usepackage{amsmath,amssymb,theorem}
\usepackage{latexsym}
\usepackage{graphicx}
\usepackage{amscd}
\usepackage{color}
\input xy
\xyoption{all}

\newcommand{\Z}{{\mathbb Z}}
\newcommand{\C}{{\mathbb{C}}}

\newcommand{\bea}{\begin{eqnarray}}
\newcommand{\eea}{\end{eqnarray}}

\newcommand{\nn}{\nonumber}

\newcommand{\bp}{\begin{pmatrix}}
\newcommand{\ep}{\end{pmatrix}}

\newcommand{\bps}{\begin{smallmatrix}}
\newcommand{\eps}{\end{smallmatrix}}

\newcommand{\la}{\langle}
\newcommand{\ra}{\rangle}

\def\s{\uparrow}

\def \cA{{\cal A}}

\def \cC{{\cal C}}

\def \cH{{\cal H}}
\def \cL{{\cal L}}
\def \cM{{\cal M}}
\def \cN{{\cal N}}
\def \cMC{{\cal MC}}

\def \cOC{{\cal OC}}

\def \V{{\cal V}}

\def \f{{\frak f}}
\def \g{{\frak g}}
\def \l{{\frak l}}
\def \m{{\frak m}}
\def \n{{\frak n}}

\def \S{{\frak S}}

\def \cb{{\bar c}}
\def \ob{{\bar o}}

\def \raw{\rightarrow}

\def \lglraw{\longleftrightarrow}

\def \oraw{\overrightarrow}

\def \End{\mathrm{End}}
\def \Hom{\mathrm{Hom}}
\def \Der{\mathrm{Der}}
\def \Coder{\mathrm{Coder}}

\def \deg{\mathrm{deg}}

\def\ott{\otimes}

\def \half{\frac{1}{2}}
\def \ov#1{\frac{1}{#1}}

\def\b1{{\mathbb 1}}

\def \fd#1{\frac{d}{d{#1}}}

\def \0{{\bf 0}}
\def \1{{\bf 1}}

\def \tri{\triangle}

\def \-{-\hspace*{-0.2cm}-}

\def\ss{\sigma}
\def \Nsddata#1#2#3#4#5{
(\xymatrix{{#1}\  \ar@<0.5ex>[r]^{{#2}} & \ {#4}
\ar@<0.5ex>[l]^{{#3}}} ,#5) }

{\theorembodyfont{\rmfamily}
 \newtheorem{defn}{Definition}
 \newtheorem{thm}{Theorem}
 \newtheorem{lem}{Lemma}

 \newtheorem{cor}{Corollary}
 \newtheorem{rem}{Remark}
}

\numberwithin{equation}{section}

\begin{document}

\begin{titlepage}
\thispagestyle{empty}
\begin{flushleft}
\hfill YITP-04-61 \\
\hfill October, 2004 \\
\hfill version 2: February, 2005\\
\hfill version 3: November, 2005\\
\end{flushleft}

\vskip 1.5 cm

\begin{center}
{\Large \textbf{Homotopy algebras inspired by}}

\vspace*{0.5cm}

{\Large \textbf{classical open-closed string field theory}}

\renewcommand{\thefootnote}{\fnsymbol{footnote}}

\vskip 2cm

{\large

\noindent{ \bigskip }\\

\it
Hiroshige Kajiura${}^{*}$ and Jim Stasheff${}^{\dag}$\\

\noindent{\smallskip  }\\

${}^{*}$
Yukawa Institute for Theoretical Physics, Kyoto University \\
Kyoto 606-8502, Japan \\
e-mail: kajiura@yukawa.kyoto-u.ac.jp\\

\noindent{\smallskip  }

${}^{\dag}$
Department of Mathematics, University of Pennsylvania \\
Philadelphia, PA 19104-6395, USA \\
e-mail: jds@math.upenn.edu\\

}
\end{center}
 
\begin{abstract}

We define a homotopy algebra associated to classical open-closed
strings. We call it an {\em  open-closed homotopy algebra (OCHA).} 
It is inspired by Zwiebach's open-closed string field theory and 
also is related to the situation of Kontsevich's deformation
quantization. 
We show that it is actually a homotopy invariant notion; 
for instance, the minimal model theorem holds. 
Also, we show that our open-closed homotopy algebra gives us 
a general scheme for deformation of 
open string structures ($A_\infty$-algebras) 
by closed strings ($L_\infty$-algebras). 
\end{abstract}

\renewcommand{\thefootnote}{}
\footnote{H.~K is supported by JSPS Research Fellowships 
for Young Scientists. 
J.~S. is supported in part by NSF grant FRG DMS-0139799 
and US-Czech Republic grant INT-0203119. 
}

\setcounter{footnote}{0}

\vfill

\end{titlepage}
\bigskip

\tableofcontents%

\section{Introduction}
\label{sec:1}

In this paper we define a strong homotopy algebra
inspired by Zwiebach's classical open-closed string field theory
\cite{Z2}
and examine its homotopy algebraic structures.
It is known that
classical closed string field theory has an
$L_\infty$-structure \cite{Z1,Sta1993,KSV} and 
classical open string field theory
has an $A_\infty$-structure  \cite{GZ,Z2,nakatsu,Ka}. 
As described by Zwiebach \cite{Z1, Z2} and others, 
string field theory is presented 
in terms of decompositions of moduli spaces of the corresponding 
Riemann surfaces into cells.
The associated Riemann surfaces are (respectively)
spheres with (closed string) punctures and
disks with (open string) punctures on the boundaries.
That is, classical closed string field theory is related to
the conformal plane $\C$ with punctures
and classical open string field theory is
related to the upper half plane $H$ with punctures on the boundary
from the viewpoint of conformal field theory.
The algebraic structure that
the classical open-closed string field theory has
is similarly interesting since it is related to
the upper half plane $H$ with punctures both in the bulk and on the
boundary,
which also appeared recently
in the context of deformation quantization \cite{Ko1,CF}.

In operad theory (see \cite{MSS}),
the relevance of the little disk operad to closed string
theory is known, where a (little) disk is related to a 
closed string puncture on a sphere in the Riemann surface picture above. 
The homology of the little disk operad
defines a Gerstenhaber algebra 
\cite{cohen, getzler-jones}, 
in particular, a suitably compatible graded commutative algebra
structure and graded Lie algebra structure. 
The framed little disk operad is in addition equipped with 
a BV-operator which rotates the disk boundary $S^1$. 
The algebraic structure on the homology is then a 
BV-algebra \cite{getz-BV}, where 
the graded commutative product and the graded Lie bracket 
are related by the BV-operator. 
Physically, closed string states associating to each disk boundary $S^1$ 
are constrained to be the $S^1$-invariant parts, 
the kernel of the BV-operator. 
This in turn leads to concentrating on the Lie algebra structure, 
where two disk boundaries are identified by {\em twist}-sewing 
as Zwiebach did \cite{Z1}. 
On the other hand, he worked at the chain level
(`off shell'), discovering an $L_\infty$-structure. This was important 
since the multi-variable operations of the $L_\infty$-structure provided
correlators of closed string field theory.
Similarly for open string theory, the little interval operad 
and associahedra are relevant, the homology corresponding 
to a graded associative algebra,
but the chain level reveals an $A_\infty$-structure 
giving the higher order correlators of open string field theory.

The corresponding operad for the open-closed string theory is the
Swiss-cheese operad \cite{Vo} that combines the little disk
operad with the little interval operad; it was inspired also by
Kontsevich's approach to deformation quantization. 
The algebraic structure at the homology level has been analyzed 
thoroughly by Harrelson \cite{erich}.
In contrast, our work in the open-closed case is 
at the level of strong homotopy algebra, 
combining the known but separate $L_\infty$- and $A_\infty$-structures.
There are interesting relations (not yet fully explored) 
between an algebra over the Swiss-cheese operad and 
the homotopy algebra we define in the present paper. 
In particular, we leave for later work the inclusion of
the appropriate homotopy algebra 
corresponding to the graded commutative product and the BV-operator.
For possible structures to be added to our structure, see \cite{Sul:oc}.

We call our structure an open-closed homotopy algebra (OCHA) 
(since it captures a lot of the operations 
in existing open-closed string field theory algebra structure 
\cite{hk-jds:physics}). 
We show that this description is a homotopy invariant algebraic
structure, i.e. 
that it transfers well under homotopy equivalences or quasi-isomorphisms. 
Also, we show that an open-closed homotopy algebra gives us
a general scheme for deformation of
open string structures ($A_\infty$-algebras)
by closed strings ($L_\infty$-algebras).

We first present our notion of 
open-closed homotopy algebras in section \ref{sec:oc}.
An open-closed homotopy algebra consists of a direct sum of
graded vector spaces $\cH =\cH_c\oplus\cH_o$. It has
an $L_\infty$-structure on $\cH_c$ and reduces to
an $A_\infty$-algebra if we set $\cH_c=0$. Moreover, the
operations that intertwine the two are a generalization of 
the strong homotopy analog of H. Cartan's notion 
of a Lie algebra $\mathfrak g$ 
acting on a differential graded algebra E \cite{cartan:notions, fgv}. 

We present the basics of these notions of homotopy algebra 
from three points of view:
multi-variable operations, coderivation differentials and tree diagrams.
In a more physically oriented paper \cite{hk-jds:physics}, 
we give an alternate interpretation in the language 
of homological vector fields on a supermanifold.

The motivating physics of string interactions suggests 
that the homotopy algebra should be appropriately {\em cyclic} 
\cite{MSS,GeK1}. 
In section 3, we make it precise 
in terms of an odd symplectic/cyclic structure which is strictly invariant 
with respect to the OCHA structure. 
It would be worth investigating a strongly homotopy invariant analog 
in the sense of \cite{tradler:infalg, tradler:infBV}, 
which we do not discuss in this paper.

One of the key theorems in homotopy algebra is the minimal model
theorem which was first proved for $A_\infty$-algebras by
Kadeishvili \cite{kadei1}. 
The minimal model theorem states the existence of minimal models 
for homotopy algebras 
analogous to Sullivan's minimal models \cite{Sul} 
for differential  graded commutative algebras 
introduced in the context of rational homotopy theory.
For an $A_\infty$- or $L_\infty$-algebra, 
the minimal model theorem is now 
combined with various stronger results; 
those employing the techniques of 
homological perturbation theory (HPT) 
(for instance see \cite{gugen,hueb-kadei,GS,gugen-lambe,GLS:chen,GLS}), 
what is called the decomposition theorem 
in \cite{thesis, KaTe}, Lef\`evre's approach \cite{Le-Ha}, etc. 
These theorems are very powerful and make clear 
the homotopy invariant nature of the algebraic properties 
(for instance \cite{markl:haha,jh-jds,KaTe}). 
In section \ref{sec:MMth} we describe these theorems 
for our open-closed homotopy algebras, pointing out subtleties 
of the open-closed case in addition to those for $L_\infty$-algebras
in comparison to the existing versions for $A_\infty$-algebras.

In section \ref{sec:deform}, we show that 
an open-closed homotopy algebra gives a general scheme of
deformation of the $A_\infty$-algebra $\cH_o$ as controlled by $\cH_c$.
A particular example of the deformation point of view 
applied in an open-closed setting occurs 
in analyzing Kontsevich's deformation quantization theorem, 
which we shall explain explicitly 
in the sequel to this paper \cite{hk-jds:physics}. 
We discuss this deformation theory also 
from the viewpoint of generalized 
Maurer-Cartan equations for an open-closed homotopy algebra and the 
moduli space of their solution space. 

We include an appendix by M.~Markl, 
where $A_\infty$-algebras over $L_\infty$-algebras are 
interpreted as a colored version of strongly homotopy algebras 
in the sense in \cite{MSS}.

We have taken care to provide the detailed signs which are crucial in
calculations, but which are conceptually unimportant and can be
ignored at first reading.
The majority of this paper is entirely mathematics, 
and in the sequel \cite{hk-jds:physics} 
we show how our structures are related to those in Zwiebach \cite{Z2}, 
deformation quantization by Kontsevich \cite{Ko1}, 
as well as those discussed in \cite{Hof1,Hof2} 
where an open-closed homotopy algebra is applied
to topological open-closed strings. 
It should be very interesting to investigate the
application to homological mirror symmetry \cite{W2,mirror,BK,Hof2}.

\section{Strong homotopy algebra}
\label{sec:oc}

An open-closed homotopy algebra, as we propose in this paper, is a strong 
homotopy algebra (or $\infty$-algebra) which combines two typical 
strong homotopy algebras, 
an $A_\infty$-algebra and an $L_\infty$-algebra. 
Let us begin by recalling those definitions. 
We restrict our arguments to the case that the characteristic of 
the field $k$ is zero. We further let $k=\C$ for simplicity. 

There are various equivalent way of defining/describing strong 
homotopy algebras: in terms of multi-variable operations and 
relations among them, in terms of a coderivation differential of 
square zero on an associated coalgebra or as a representation of a 
particular operad of trees. 
We will treat all three of these in turn. 
The reader who is familiar with these approaches 
to the `classical' $A_\infty$-algebras and  $L_\infty$-algebras 
can move ahead to subsection \ref{ssec:ocha}, 
being warned that the definitions 
we give are different from the original ones \cite{Sta,LS,MSS} 
in the degrees of the multi-linear maps and hence of the relevant signs. 
Both are in fact equivalent and related by  
{\it suspension} \cite{MSS}, as we explain further below.

\subsection{Strong homotopy associative algebras}
\begin{defn}[$A_\infty$-algebra\ 
(strong homotopy associative algebra)\cite{Sta}]
Let $A$ be a $\Z$-graded vector space 
$A=\oplus_{r\in\Z} A^r$ and suppose that 
there exists a collection of degree one multi-linear maps 
\begin{equation*}
 \m:=\{m_k : A^{\otimes k}\raw A\}_{k\ge 1} \ .
\end{equation*}
$(A,\m)$ is called an {\em  $A_\infty$-algebra} when the multi-linear
maps $m_k$ satisfy the following relations 
\begin{equation}
\sum_{k+l=n+1}\sum_{i=1}^{k}
{(-1)^{o_1+\cdots+o_{i-1}}
 m_k(o_1,\cdots,o_{i-1},m_l(o_{i},\cdots,o_{i+l-1}),
 o_{i+l},\cdots,o_n)}=0\ 
 \label{Ainfty}
\end{equation}
for $n\ge 1$, 
where $o_j$ on $(-1)$ denotes the degree of $o_j$. 

A {\em weak  $A_\infty$-algebra} consists of a 
collection of degree one multi-linear maps 
\begin{equation*}
 \m:=\{m_k : A^{\otimes k}\raw A\}_{k\ge 0} 
\end{equation*} 
satisfying the above relations, but for $n\ge 0$ and 
in particular with $k,l\ge 0$.
 \label{defn:Ainfty}
\end{defn}
\begin{rem}
The relation (\ref{Ainfty}) is different from the original one \cite{Sta}
in the definition of the degrees of the multi-linear maps $m_k$
and hence of the signs.
Both are in fact equivalent and related by
{\it desuspension} \cite{MSS}.
In \cite{Sta}, the $m_k$ are multi-linear maps on ${\downarrow}A$ 
where $({\downarrow}A)^{r+1} = A^r$; we denote desuspension by
$\downarrow$. (The algebraic geometry tradition would use $[-1]$. )
Note that, in that notation \cite{Sta}, a differential graded (dg) 
algebra is an $A_\infty$-algebra with a differential $m_1$, 
a product $m_2$, and $m_3=m_4=\cdots=0$.

The `weak' version is fairly new, inspired by physics, where 
$m_0:\C\to A,$  regarded 
as an element $m_0(1)\in A$, is related to what physicists 
refer to  as a `background'.
The augmented relation then implies that 
$m_0(1)$ is a cycle, but $m_1 m_1$ need no longer be 0,  rather
$m_1 m_1 = \pm m_2(m_0\otimes 1) \pm m_2(1\otimes m_0)$. 
 \label{rem:sus}
\end{rem}
\begin{defn}[$A_\infty$-morphism]
For two $A_\infty$-algebras $(A, \m)$ and $(A',\m')$,
suppose that there exists a collection of
degree zero (degree preserving) multi-linear maps
\begin{equation*}
 f_k: A^{\otimes k}\raw A'\ ,\qquad k\geq 1\ .
\end{equation*}
The collection
$\{f_k\}_{k\ge 1}:(A,\m)\raw (A',\m')$ is called an {\em
$A_\infty$-morphism} iff it satisfies the following relations:
\begin{equation}
 \begin{split}
& \sum_{1\leq k_1<k_2\dots <k_j=n}\
{m'_j(f_{k_1}(o_1,\cdots,o_{k_1}),
f_{k_2-k_1}(o_{k_1+1},\cdots,o_{k_2}),\cdots, f_{n-k_{j-1}}
(o_{k_{j-1}+1},\cdots,o_n))}\\
&\qquad=\sum_{k+l=n+1}\sum_{i=1}^{k}{(-1)^{o_1+\dots +o_{i-1}}
f_k(o_1,\cdots,o_{i-1},m_l(o_{i},\cdots,o_{i+l-1}),
o_{i+l},\cdots,o_n)}\
\end{split}
\label{amorphism}
\end{equation}
for $n\ge 1$.

If $(A, \m)$ and $(A',\m')$ are {\em weak} $A_\infty$-algebras, then a
{\em weak} $A_\infty$-morphism consists of multi-linear maps
$\{f_k\}_{k\ge 0}$, where $f_0:\C\raw A'$, satisfying 
the above conditions 
for $n\ge 0$. In particular, the condition for $n=0$ is: 
$$
 f_1\circ m_0 = \sum_{k\ge 0} m'_k(f_0,\cdots,f_0)\ .
$$
 \label{defn:Ainftymorp}
\end{defn}

 \subsection{The coalgebra description and the Gerstenhaber bracket}
\label{ssec:gerst}

The maps $m_k$ can be assembled into a single map, also
denoted $\m$, from the tensor space $T^cA = \oplus_{k\geq 0}
A^{\ott k}$ to
$A$  with the convention that  $A^{\ott 0}=\C$. The grading implied
by having the maps $m_k$ all of degree one is the usual grading
on each $A^{\ott k}.$ We can regard
$T^cA$ as the tensor {\em coalgebra} by defining
\begin{equation*}
 \tri(o_1 \ott\cdots\ott o_n) 
 = \Sigma_{p=0}^n (o_1\ott\cdots\ott o_p)\ott (o_{p+1}\ott\cdots\ott o_n)\ .
\end{equation*}
A map $f\in \Hom(T^cA, T^cA)$
is a {\em graded coderivation} means 
$\tri f = (f\otimes \1 + \1\otimes f)\tri$, 
with the appropriate signs and dual to the definition of a
graded derivation of an algebra. 
Here $\1$ denotes the identity $\1:A\to A$. 
We then identify
$\Hom(T^cA, A)$ with $\Coder(T^cA)$ by lifting a multi-linear map
as a coderivation \cite{jds:intrinsic}.  Analogously to the situation for
derivations, the composition graded  commutator of coderivations is
again a coderivation; this graded  commutator corresponds 
to the {\em Gerstenhaber bracket} on 
$\Hom(T^cA, A)$ \cite{gerst:coh,jds:intrinsic}. 
Notice that this involves a shift in grading 
since Gerstenhaber uses the traditional Hochschild complex grading. 
Thus $\Coder(T^cA)$ is a graded Lie algebra 
and in fact a dg Lie algebra 
with respect to the bar construction differential, which corresponds 
to the Hochschild differential on $\Hom(T^cA, A)$ 
in the case of an associative algebra $(A,m)$\cite{gerst:coh}. 
Using the bracket, the differential can be written as $[m,\ ]$.

The advantage of this point of view is that the component maps $m_k$
assemble into a single map $\m$ in $\Coder(T^cA)$ 
and relations (\ref{Ainfty}) can be summarized by 
\begin{equation*}
 [\m,\m]=0\  {\rm\quad or, equivalently, \quad} D^2=0\ , 
\end{equation*}
where $D=[\m, \ \,]$. 
In fact, $\m\in \Coder(T^c A)$ is an $A_\infty$-algebra structure 
on $A$ iff $[\m,\m]=0$
{\em and} $\m$ has no constant term, $m_0 =0$. 
If $m_0\neq 0$, the structure is a weak $A_\infty$-algebra.
The $A_\infty$-morphism components similarly combine to give a single map
of dg coalgebras $\f :T^cA \raw T^cA'$, 
$(\f\otimes\f)\tri=\tri\f$. 
In particular, (\ref{amorphism})
is equivalent to $\f \circ\m = \m'\circ \f$.

 \subsection{The tree description}
\label{ssec:Ainftytree}

\begin{figure}[h]
\begin{equation}
\begin{minipage}[c]{50mm}{\includegraphics{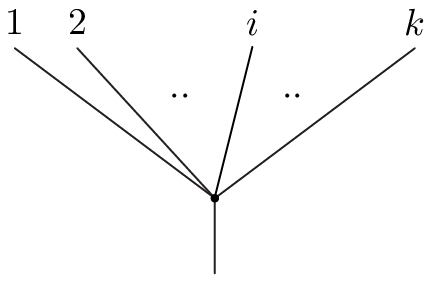}}
\end{minipage}
\bullet_i
\begin{minipage}[c]{40mm}{\includegraphics{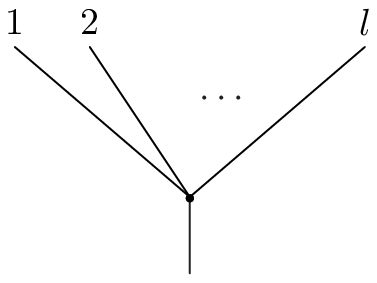}}
\end{minipage}\ =\ 
\begin{minipage}[c]{50mm}{\includegraphics{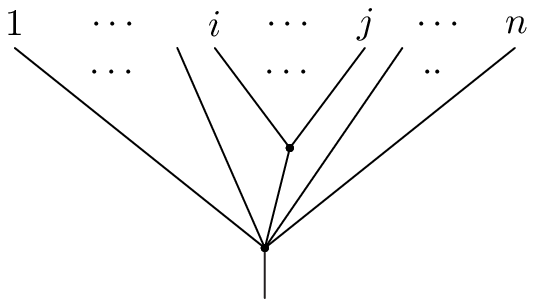}}
\end{minipage}
 \label{m-graft}
\end{equation}
 \caption{The grafting $m_k\bullet_i m_l$ of the $l$-corolla 
$m_l$ to the $i$-th leaf of $k$-corolla $m_k$, 
where $j=i+l-1$ and $n=k+l-1$. }
 \label{fig:grafting}
\end{figure}

There are some advantages to indexing the maps $m_k$ and 
their compositions by planar rooted trees; 
e.g. $m_k$ will correspond to the {\em corolla} 
with $k$ leaves all attached directly to the root. 
The composite $m_k\bullet_i m_l$ then corresponds to grafting the root of
$m_l$ to the $i$-th leaf of $m_k$, reading from left to right 
(see Figure \ref{fig:grafting}). 
This is the essence of the planar rooted tree operad \cite{MSS}. 
Multi-linear maps compose in just this way, 
so relations (\ref{Ainfty}) can be phrased 
as saying we have a map from planar rooted trees
to multi-linear maps respecting the $\bullet_i$ `products', 
the essence of a map of operads \cite{MSS}. 
More precisely, let $\cA_\infty(n)$, $n\ge 1$ 
be a graded vector space 
spanned by planar rooted trees of $n$ leaves 
with identity $e\in \cA_\infty(1)$. 
For a planar rooted tree $T\in\cA_\infty(n)$, 
its grading is introduced by 
the number of the vertices contained in $T$, 
which we denote by $v(T)$. 
A tree $T\in \cA_\infty(n)$, $n\ge 2$, with $v(T)=1$ is 
the corolla $m_n$. 
Any tree $T$ with $v(T)=2$ is obtained by the grafting of 
two corollas as in eq.(\ref{m-graft}). 
Grafting of any two trees is defined in a similar way, 
with an appropriate sign rule, 
and any tree $T$ with $v(T)\ge 2$ can be obtained recursively by 
grafting a corolla to a tree $T'$ with $v(T')=v(T)-1$. 
One can define a differential $d$ of degree one, which acts on 
each corolla as 
\begin{equation}
 d\left( m_n \right)=- \sum_{k,l\ge 2,\ k+l=n+1}\sum_{i=1}^{k}
 m_k\bullet_i m_l
 \label{Ainfty-operad}
\end{equation}
and extends to one on $\cA_\infty:=\oplus_{n\ge 1}\cA_\infty(n)$ by 
the following rule: 
\begin{equation*}
 d(T\bullet_i T')= d(T)\bullet_i T' +(-1)^{v(T)} T\bullet_i d(T')\ . 
\end{equation*}
If we introduce the contraction of internal edges, that is, 
indicate by $T'\raw T$ that $T$ is obtained from $T'$ by 
contracting an internal edge, the differential is equivalently given 
by 
\begin{equation*}
 d(T)= \sum_{T'\raw T}\pm T'\ 
\end{equation*}
with an appropriate sign $\pm$. 
Thus, one obtains a dg operad $\cA_\infty$, which is known as 
the {\em $A_\infty$-operad}. 
An algebra $A$ over $\cA_\infty$ is obtained by 
a representation $\phi:\cA_\infty(k)\raw \Hom(A^{\otimes k},A)$, 
though we use the same notation $m_k$ 
for the $k$-corolla and its image by $\phi$. 
This map $\phi$ should be defined 
so that it is compatible with respect to $\bullet_i$ and 
also the differentials. 
Here, identifying elements in $\Hom(A^{\otimes k},A)$ 
with those in $\Coder(T^cA)$, 
the differential on the algebra side is $[m_1,\ ]$. 
In particular, for each corolla one gets 
\begin{equation*}
 \phi \left(d(m_k)\right)= m_1\phi(m_k)+ 
\sum_{i=1}^k \phi(m_k)\circ
 (\1^{\otimes (i-1)}\otimes m_1\otimes \1^{\otimes (k-i)})\ .
\end{equation*}
It is clear that this equation combined with eq.(\ref{Ainfty-operad}) 
implies the $A_\infty$-relations (\ref{Ainfty}).  

Note that this grading of trees we introduced here 
corresponds to that in the suspended notation 
of $A_\infty$-algebras. In unsuspended notation ${\downarrow}A$, 
the grading of a tree $T\in\cA_\infty(n)$ should be replaced by 
$int(T)+2-n$, where $int(T)=v(T)-1$ denotes the number of the 
internal edges. 
Each tree $T\in\cA_\infty(n)$ corresponds 
to a codimension $int(T)$ boundary piece of associahedron $K_n$ \cite{Sta}, 
that is, $int(T)+2-n$ is equal to minus the dimension of the 
boundary piece.

\subsection{$L_\infty$-algebras and morphisms}

Since an ordinary Lie algebra $\g$ is regarded as ungraded, the
defining  bracket is regarded as skew-symmetric. 
If we regard $\g$ as all of degree one,
then the bracket would be graded symmetric.
For dg Lie algebras and $L_\infty$-algebras, we need graded symmetry,
which refers to symmetry with signs determined by the grading. 
The basic relation is 
\begin{equation}
\tau :x\otimes y \mapsto (-1)^{|x||y|}y\ott x\ .
\label{tau}
\end{equation}
Also we adopt the convention that tensor products of graded functions or
operators
have the signs built in; e.g.
$(f\ott g)(x\ott y) = (-1)^{|g| |x|}f(x)\ott g(y).$
By decomposing permutations as a product of transpositions, there is 
then defined the sign of a permutation of $n$ graded elements, e.g
for any $c_i\in V$ $1\le i\le n$ and any $\sigma\in\S_n$, 
the permutation of $n$ graded elements, is defined by
\begin{equation}
 \sigma(c_1,\cdot\cdots,c_n)
 =(-1)^{\epsilon(\sigma)}(c_{\sigma(1)},\cdots\cdot,c_{\sigma(n)})\ .
\end{equation}
The sign  $(-1)^{\epsilon(\sigma)}$ is often referred to as the Koszul
sign of the permutation.
\begin{defn}[Graded symmetry]
A {\em graded symmetric multi-linear map}  of a
graded vector space $V$ to itself is a linear map
$f:V^{\otimes n}\to V$ such that for any $c_i\in V$, $1\le i\le n$, and
any $\sigma\in\S_n$ (the permutation group of $n$ elements), the relation
\begin{equation}
 f(c_1,\cdot\cdots,c_n)
 =(-1)^{\epsilon(\sigma)}f(c_{\sigma(1)},\cdots\cdot,c_{\sigma(n)})
\end{equation}
holds.
\end{defn}
Since we will have many formulas with such indices and their
permutations, we will use the notation $I=(i_1,\dots,i_n)$ and
\begin{equation*}
 c_I = c_{i_1}\ott \cdots \ott c_{i_n}\ .
\end{equation*}
Then, for any $\sigma\in\S_n$, we use $\sigma(I)$ to denote
$(\sigma(i_1),\dots,\sigma(i_n))$
and hence
\begin{equation*}
 c_{\sigma(I)} = c_{\sigma(i_1)}\ott \cdots\ott c_{\sigma(i_n)}\ .
\end{equation*}
\begin{defn} By a $(k,l)${\em -unshuffle} of
$c_1,\dots,c_n$ with $n=k+l$ 
is meant a permutation $\sigma$ such that for $i<j\leq k$, 
we have
$\sigma(i)< \sigma(j)$ and similarly for $k < i<j\leq k+l$. 
We denote the subgroup of $(k,l)$-unshuffles in $\S_{k+l}$
by $\S_{k,l}$ and by $\S_{k+l=n},$ the union of the subgroups
$\S_{k,l}$ with $k+l=n$. Similarly, a $(k_1,\cdots,k_i)$-unshuffle
means a permutation $\sigma\in \S_n$ with $n=k_1+\cdots+k_i$ such that
the order is preserved within each block of length $k_1,\cdots,k_i$. 
The subgroup of $\S_n$ consisting of all such unshuffles 
we denote by $\S_{k_1,\cdots,k_i}$. 
\end{defn}
\begin{defn}[$L_\infty$-algebra\ (strong homotopy Lie algebra)
\cite{LS}]
Let $L$ be a graded vector space and suppose that
a collection of degree one graded symmetric linear maps
$\l:=\{l_k:L^{\otimes k}\raw L\}_{k\ge 1}$ is given. 
$(L,\l)$ is called an {\em  $L_\infty$-algebra} iff
the  maps satisfy the following relations:
\begin{equation}
\sum_{\sigma\in\S_{k+l=n}}
(-1)^{\epsilon(\sigma)}
l_{1+l}(l_k(c_{\sigma(1)},\cdots,c_{\sigma(k)}),
 c_{\sigma(k+1)},\cdots,c_{\sigma(n)})=0\
 \label{Linfty}
\end{equation}
for $n\ge 1$. Using the multi-index notation $I$, this can be written
\begin{equation}
\sum_{\sigma\in\S_{k+l=n}}
(-1)^{\epsilon(\sigma)}
l_{1+l}(l_k\otimes \1^{\otimes l})(c_{\sigma(I)} )=0\
 \label{LinftyI}
\end{equation}
for $n\geq 1$. 
A {\em weak $L_\infty$-algebra} consists of 
a collection of degree one graded symmetric linear maps
$\l:=\{l_k: L^{\otimes k}\raw L\}_{l\ge 0}$
satisfying the above relations, but for $n\ge 0$ and 
with $k,l\ge 0$.
 \label{defn:Linfty}
\end{defn}
\begin{rem}
The alternate definition in which the summation is 
over all permutations, rather than just unshuffles, requires 
the inclusion of appropriate coefficients involving factorials. 
Recall that the signs we use correspond to the 
suspension of the original definition. 
\end{rem}
\begin{rem}
A dg Lie algebra is expressed as the desuspension of an 
$L_\infty$-algebra $(L,\l)$ where 
$l_1$ and $l_2$ correspond to the differential and the Lie bracket, 
respectively, and higher multi-linear maps $l_3,l_4,\cdots$ are absent. 
 \label{rem:susL}
\end{rem}
\begin{rem}
For the `weak' version, remarks analogous to those 
for weak $A_\infty$-algebras apply.
 \label{rem:weakL}
\end{rem}
\begin{defn}[$L_\infty$-morphism]
For two $L_\infty$-algebras $(L, \l)$ and $(L',\l')$,
suppose that there exists a collection of
degree zero (degree preserving)
graded symmetric multi-linear maps
\begin{equation*}
 f_k: L^{\otimes k}\raw L'\ ,\qquad l\ge 0\ .
\end{equation*}
Here $f_0$ is a map from $\C$ to a degree zero subvector space of
$L$. The collection
$\{f_k\}_{k\ge 1}:(L,\l)\raw (L',\l')$ is called an 
{\em $L_\infty$-morphism} iff it satisfies the following relations
\begin{equation}
 \begin{split}
 \sum_{\sigma\in\S_{k+l=n}} &(-1)^{\epsilon(\sigma)}
f_{1+l}(l_k\otimes \1^{\otimes l})(c_{\sigma(I)}) \\
&=\sum_{\sigma\in\S_{k_1+\cdots +k_j=n}}
\frac{(-1)^{\epsilon(\sigma)}}{j!}
l'_j(f_{k_1}\otimes f_{k_2}\otimes \cdots \otimes
f_{k_j})(c_{\sigma(I)})
\end{split}
\label{Linftymorp}
\end{equation}
for $n\ge 1$. 

When $(L, \l)$ and $(L',\l')$ are {\em weak} $L_\infty$-algebras, 
then a {\em weak} $L_\infty$-morphism consists of multi-linear maps
$\{f_k\}_{k\ge 0}$ satisfying the above conditions
and in addition $f_1\circ l_0 = \sum_k (1/k!) l'_k(f_0,\cdots,f_0)$. 
 \label{defn:Linftymorp}
\end{defn}
When $L'$ is (a suspension of) a strict dg Lie algebra, 
the formula simplifies greatly
since, on the right hand side, we have $j=1$ or $2$ only
(\cite{lada-markl} Definition 5.2).

\subsection{The symmetric coalgebra description}
\label{ssec:symm}

The graded symmetric coalgebra on a graded vector space $V$ is naturally
the {\em sub}coalgebra $S^c V\subset T^c V$ 
consisting of the graded symmetric elements in each $V^{\ott n}$. 
By {\em not} regarding $S^c V$ as a quotient of $T^c V$, 
certain complicated factorial coefficients do not
appear in our formulas. Also, in rational homotopy theory, $S^c V$ is
often denoted $\Lambda^c V$, due to a historical accident. To avoid
possible confusion, we will use neither, but instead $C(V)$, 
as in \cite{SS}.

Again, the sum of the maps $\l =\oplus_k l_k$ provides a coderivation
differential $\l$ (with $[\l,\l]=0$) on the full tensor coalgebra.
Because of the graded symmetry of the $l_k$, 
the structure can be identified with a coderivation differential on the
graded symmetric coalgebra $C(L)$, see \cite{lada-markl}. That is,
$\l\in \Coder^1(C(L))$ is an $L_\infty$-algebra structure on $L$ iff
$[\l,\l]=0$
{\em and} $\l$ has no constant term: $l_0 =0$. If 
$l_0\neq 0$, the structure is a weak $L_\infty$-algebra.

Also, 
if we assemble the $L_\infty$-morphism components $\{f_k\}$ 
to a single $\f:T^c L\to L'$ and lift it
to a {\em coalgebra} map $\f:C(L)\to C(L')$, 
then (\ref{Linftymorp}) is equivalent to $\f\circ\l = \l'\circ\f$.

 \subsection{The tree description}
\label{ssec:Linftytree}

\begin{figure}[h]
\begin{equation}
\begin{minipage}[c]{50mm}{\includegraphics{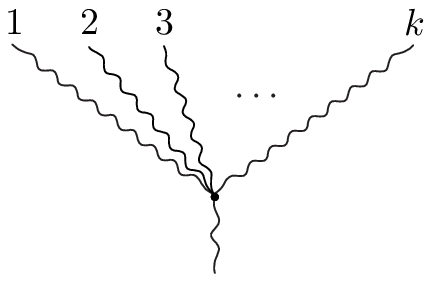}}
\end{minipage}
=
\begin{minipage}[c]{50mm}{\includegraphics{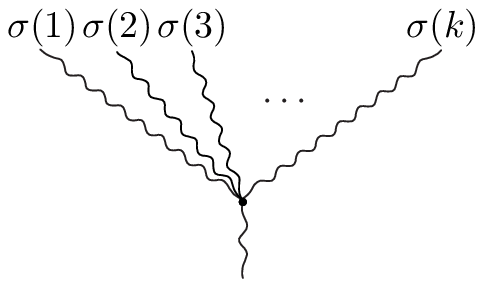}}
\end{minipage}
\end{equation}
 \caption{Nonplanar $k$-corolla corresponding to $l_k$. 
Since edges are non-planar, it is symmetric with respect to 
the permutation of the edges. }
 \label{fig:l}
\end{figure}

The tree operad description of $L_\infty$-algebras uses
non-planar rooted trees with leaves numbered $1,2,...$ arbitrarily 
\cite{MSS}. 
Namely, 
a non-planar rooted tree can be expressed as a planar rooted tree 
but with arbitrary ordered labels for the leaves. 
In particular, corollas obtained by permuting the labels are 
identified (Figure \ref{fig:l}).  
Let $\cL_\infty(n)$, $n\ge 1$ be a graded vector space generated by 
those non-planar rooted trees of $n$ leaves. 
For a tree $T\in\cL_\infty(n)$, a permutation $\sigma\in\S_n$ 
of the labels for leaves generates a different tree in general, 
but sometimes the same one because of the symmetry 
of the corollas above. 
The grafting, $\circ_i$, to the $i$-th leaf is defined  as in the 
planar case in subsection \ref{ssec:Ainftytree}, 
and any non-planar rooted tree is obtained by 
grafting corollas $\{l_k\}_{k\ge 2}$ recursively, as in the planar case, 
together with the permutations of the labels for the leaves. 
A degree one differential 
$d: \cL_\infty(n)\raw\cL_\infty(n)$ is given in a similar way; 
for $T'\raw T$ indicating that $T$ is obtained from $T'$ 
by the contraction of an internal edge, 
\begin{equation*}
 d(T)=\sum_{T'\raw T}\pm T'\ .
\end{equation*}
In particular, for each corolla one gets 
\begin{equation}
d \left(\, 
\begin{minipage}[c]{45mm}{\includegraphics{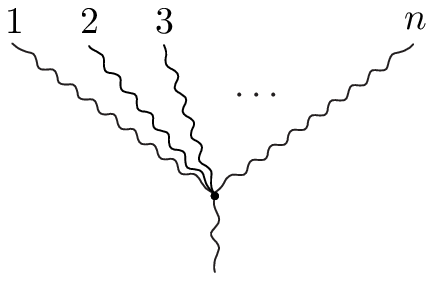}}
\end{minipage}
\ \right)
\ = - \ \sum_{\sigma\in\S_{k+l=n}}
\begin{minipage}[c]{50mm}{\includegraphics{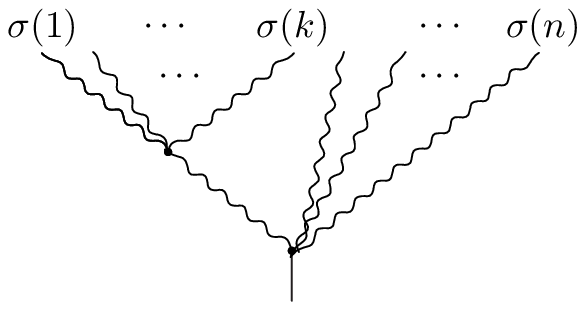}}
\end{minipage}\ ,
 \label{d_Linfty-operad}
\end{equation}
and $d(T\circ_i T')= d(T)\circ_i T' +(-1)^{v(T)} T\circ_i d(T')$ 
again holds. 
Thus, $\cL_\infty:=\oplus_{n\ge 1}\cL_\infty(n)$ forms a dg operad, 
called the {\em $L_\infty$-operad}. 

An algebra $L$ over $\cL_\infty$ obtained by a map 
$\phi:\cL_\infty(k)\raw \Hom(L^{\otimes k},L)$ then forms an 
$L_\infty$-algebra $(L,\l)$. 
If we use double desuspended notation ${\downarrow}{\downarrow}L$ 
(physics notation \cite{LS,Z1}; see \cite{KSV}), 
the degree of a multi-linear map $l_k$ turns into  
$1-2(k-1)=3-2k$. 
The grading of a tree $T\in\cL_\infty(n)$ should be replaced by 
$int(T)+(3-2k)$, where 
the dimension of the 
corresponding boundary piece of the compactified moduli space 
of a sphere with $(k+1)$ marked points is $-(int(T)+(3-2k))-1$ 
(see \cite{KSV}).

 \subsection{Open-closed homotopy algebra (OCHA)}
\label{ssec:ocha}

For our open-closed homotopy algebra, we consider a graded vector space
$\cH=\cH_c\oplus\cH_o$ in which $\cH_c$ will be an $L_\infty$-algebra 
and $\cH_o$, an $A_\infty$-algebra. 
An open-closed homotopy algebra includes various sub-structures, 
or reduces to various simpler structures as particular cases. 
An important such structure is the action of
$\cH_c$ as an $L_\infty$-algebra on $\cH_o$ as a dg vector space. 
This is the appropriate strong homotopy version of the action 
of an ordinary Lie algebra $L$ on a vector space $M$, also described 
as $M$ being a module over $L$ or a representation of $L$. 
Thus we can also speak of $\cH_o$ as a {\em strong homotopy module} 
over $\cH_c$ or as a {\em strong homotopy representation} of $\cH_c$
(cf. \cite{jds:bull}).
Moreover, we will need the strong
homotopy version of an algebra $A$ over a Lie algebra $L$, that is, an
action of $L$ by derivations of $A,$ so that the  map $L\to \End(A)$
takes values in the Lie sub-algebra $\Der A$.
We first arrange these known structures, 
and then define our open-closed homotopy algebra 
as an extension of them. 

Lada and Markl (\cite{lada-markl} Definition 5.1) provide 
the definition of an $L_\infty$-module 
at the desuspended level where it is easier to see the structure 
as satisfying the relations for a Lie module only up to homotopy. 
At the suspended level with which we are working in this paper, 
adjusting degrees and signs, the definition looks as follows:
\begin{defn}[sh-L-module]
\label{snuff}
Let $L=(L,l_i)$, be an $L_\infty$-algebra, and let $M$ be 
a differential graded vector space with differential denoted by $k_1$. 
Then a {\em left $L_\infty$-module structure on $M$} 
is a collection $\{k_n| 1 < n < \infty\}$ of
graded linear maps of degree one 
$$
 k_n: L^{\otimes (n-1)}\otimes M \longrightarrow M\ ,
$$
which are graded symmetric in $L$ and such that
\begin{equation}
 \begin{split}
&\sum_{\sigma\in\S_{p+q=n}}(-1)^{\epsilon(\sigma)}
k_{1+q+1}(l_p(\xi_{\sigma(1)},\cdots,\xi_{\sigma(p)}),
\xi_{\sigma(p+1)},\cdots,\xi_{\sigma(n)},\xi_{n+1}) \\
& +\sum_{\sigma\in\S_{p+q=n}}
(-1)^{\epsilon(\sigma)}
(-1)^{\sigma(1)+\cdots +\sigma(p)}
k_{p+1}(\xi_{\sigma(1)},\cdots,\xi_{\sigma(p)},
k_{q+1}(\xi_{\sigma(p+1)},\cdots,\xi_{\sigma(n)},\xi_{n+1}))=0\ ,
 \end{split}
 \label{shL}
\end{equation}
where $\xi_1,\cdots,\xi_n\in L$ and $\xi_{n+1}\in M$.
 \label{defn:sh-module}
\end{defn}
Several comments are in order.
For clarity, let $d_L=l_1$ and $d_M=k_1$, 
then the first few relations are: 
\begin{itemize}
\item $d_M k_2=-k_2 (d_L\ott\1+\1\ott d_M)$, that is, 
$k_2$ is a chain map\ ,
\item $d_M k_3+k_3 (d_L\ott\1\ott\1+\1\ott d_L\ott\1
+ \1\ott \1\ott d_M)
=-k_2(l_2\ott\1+(\1\otimes k_2)\circ(\tau\ott\1)+\1\ott k_2)$\ , 
\end{itemize}\noindent
where $\tau$ is the interchange operator (\ref{tau});
that is, $k_2$ is a Lie action up to homotopy.

Of course, the fundamental example of such a structure occurs in
the situation in which $M=L$ and each $k_i=l_i$, i.e., $L$ is an
$L_\infty$-module over itself.

According to Lada and Markl \cite{lada-markl}, we have the usual
relationship between
homomorphisms and module structures. Let $\End(M)$ denote the
differential
graded Lie algebra of linear maps from $M$ to $M$ with
bracket given by the  composition graded commutator and differential
induced by the
differential $k_1$ on $M$. After shifting the grading, their theorem
reads as follows:
\begin{thm}[Theorem 5.4 \cite{lada-markl}]
\label{corresp}
Suppose that $L$ is an $L_\infty$-algebra and that $M =(M,k_1)$ is
a differential graded vector space. Then there exists a natural
one-to-one
correspondence between sh-L-module
structures on $M$ and $L_\infty$-maps $L\to {\s}\End(M)$.
\end{thm}
\begin{rem} 
The following is phrased with the traditional grading. Fortunately,
$\End(M)$ and $\End({\s}M)$ are isomorphic as graded Lie algebras, as
are the Lie subalgebras $\Der A$ and $\Der ({\s}A)$. 
\end{rem}
In his ground breaking ``Notions d'alg\`ebre diff\'erentielle; ...''
\cite{cartan:notions}, 
Henri Cartan formalized several dg algebra notions related to his study of
the deRham cohomology of principal fibre bundles, in particular, that
of a Lie group $G$ acting in (`dans') a differential graded algebra E. 
The action uses only the Lie algebra $\mathfrak g$ of $G$. 
Cartan's action includes both graded derivations, the Lie derivative
$\theta (X)$ and the inner derivative $i(X)$ 
for $X\in \mathfrak g$. 
We need only the analog of the $\theta (X)$,
(which we denote $\rho(X)$ since by $\theta$ we denote the image 
by $\rho$ of an element $X\in\g$) 
for the following definition.
\begin{defn} For a dg Lie algebra $\mathfrak g$, 
a dg associative algebra $A$ is a $\mathfrak g${\em -algebra} 
if $\mathfrak g$ acts by derivations of $A,$ i.e. 
there is given a representation of $\mathfrak g$, 
\begin{equation*}
 \rho:\mathfrak g \to \Der A
\end{equation*}
which is a Lie map and a chain map. 
 \label{defn:g-alg}
\end{defn}
If we denote $\rho(X)a$ as $Xa$, then this means
\begin{equation}
X(ab) = (Xa)b \pm aX(b)
 \label{Xab}
\end{equation}
\begin{equation}
[X,Y]a = X(Ya)\pm Y(Xa)
 \label{XYa}
\end{equation}
\begin{equation}
d_A(Xa) = (d_{\mathfrak g}X)a \pm X(d_Aa)\ .
\end{equation}
\begin{rem} 
The concept was later reintroduced 
by Flato, Gerstenhaber and Voronov \cite{fgv} 
under the name {\em Leibniz pair}, cf. also \cite{alekseev-meinrenken}.
\end{rem}
We have seen that an $L_\infty$-module structure is defined 
in terms of relations on the maps $L^{\otimes p}\otimes A\raw A$. 
A $\mathfrak g$-algebra (or Leibniz pair) extends this
in the sense that it includes a relation (\ref{Xab}) 
on $\g\otimes A^{\otimes 2}\raw A$ 
where $\g={\downarrow}L$. 
For an extension to $L\otimes A^{\otimes q}\raw A$, 
a relevant notion is that of a homotopy derivation; 
that is, given $\theta_1:A \to A$ and $m_2: A\otimes A \to A$,
we ask for a homotopy $\theta_2:A \otimes A \to A$ between 
$\theta_1 m_2$ and $m_2(\theta_1\otimes \1) + m_2(\1\otimes\theta_1)$.  
Further higher homotopies follow the usual pattern.
\begin{defn}[Strong homotopy derivation]
A {\em strong homotopy derivation} of degree one of an $A_\infty$-algebra
$(A,\m)$ consists of a collection of multi-linear maps of degree one 
\begin{equation*}
 \theta:=\{\theta_q : A^{\otimes q}\raw A\}_{q \geq 1}
\end{equation*}
satisfying the following relations
\begin{equation}
\begin{split}
&0=\sum_{r+s=q+1}
\sum_{i=0}^{r-1}(-1)^{\beta(s,i)}\theta_r
(o_1,\cdots,o_i,m_s(o_{i+1},\cdots,o_{i+s}),\cdots,o_q)
\\
&\quad \quad\quad\hskip10ex
+ (-1)^{\beta(s,i)}
m_r(o_1,\cdots,o_i,\theta_s(o_{i+1},\cdots,o_{i+s}),\cdots,o_q)\ .
\end{split}
\label{shder}
\end{equation}
 \label{defn:shder}
\end{defn}
Here the sign $\beta(s,i)=o_1+\cdots +o_i$ results
from moving $m_s$, respectively $\theta_s,$ past $(o_1,\cdots,o_i).$
\begin{rem}
The formulas are equivalent (as suggested by Markl \cite{markl:1203})
to seeing $\theta$ as a
coderivation $\theta$ of $T^cA$ with no constant term and such that
\begin{equation*}
 [\m,\theta] = 0\ .
\end{equation*}
\end{rem}
If we extend $\theta_q: A^{\otimes q}\raw A$ to a map 
$\rho_q: L\otimes A^{\otimes q}\raw A$ by $\theta_q:=\rho_q({\s}X)$, 
the appropriate defining equation is then replaced by 
\begin{equation}
 \rho(d_\g({\s}X)) = [\m, \rho({\s}X)]\ ,
 \label{rho-chain}
\end{equation}
where $\rho({\s}X)$ is the lift of $\sum_q\rho_q({\s}X)$ to a coderivation.   
This can be read as the condition for $\rho$ to be 
a chain map regarding $[\m, \ ]$ as a 
differential on ${\s}\Coder(T^cA)$.

Now if $A$ is an $L_\infty$-module over $L$, the analog of Cartan's
second condition would be for $\rho: L \to {\s}\Coder(T^c A)$ to be 
an $L_\infty$-map. 
We already have the homotopies
$k_3:L\ott L\ott A \to A$ and
$\rho_2:L \otimes A \otimes A \to A$. The next stage of a {\em strong}
homotopy
version
looks at the various compositions giving rise to maps $L^{\ott p}\ott
A^{\ott q} \to A$ with
$p+q = 4$ and so forth.
\begin{defn}[$A_\infty$-algebra over an $L_\infty$-algebra]
Let $L$ be an $L_\infty$-algebra and $A$ an $A_\infty$-algebra 
which as a dg vector space is an sh-L module. 
That $A$ is an {\em $A_\infty$-algebra over $L$} means 
that the module structure map
$\rho:L\to {\s}\End(A)$, regarded as in ${\s}\Coder(T^cA)$, extends to
an $L_\infty$-map $L\to {\s}\Coder(T^cA)$, 
where ${\s}\Coder(T^cA)$ is the suspension of $\Coder(T^cA)$ 
as the dg Lie algebra stated in subsection \ref{ssec:gerst}. 
 \label{defn:AovL}
\end{defn}
\begin{thm}
That $A$ is an {\em $A_\infty$-algebra over the $L_\infty$-algebra $L$} is 
equivalent to having a family of maps 
$\n=\{n_{p,q}:L^{\ott p}\ott A^{\ott q} \to A\}$ for $p\geq 0$ 
but $q>0$ satisfying the compatibility conditions:
\begin{equation}
\begin{split}
&0=\sum_{\sigma\in\S_{p+r=n}}(-1)^{\epsilon(\sigma)}
n_{1+r,m}(l_p(c_{\sigma(1)},\cdots,c_{\sigma(p)}),
c_{\sigma(p+1)}\cdots,c_{\sigma(n)};o_1,\cdots,o_m) \\
& +
\sum_{\substack{\sigma\in\S_{p+r=n}\\i+s+j=m}}
(-1)^{\mu_{p,i}(\sigma)}
n_{p,i+1+j}(c_{\sigma(1)},\cdot\cdot,c_{\sigma(p)};o_1,\cdot\cdot,o_i,
n_{r,s}(c_{\sigma(p+1)},\cdot\cdot,c_{\sigma(n)};o_{i+1},\cdot\cdot,o_{i+s}),
o_{i+s+1},\cdot\cdot,o_m) \ . 
\end{split}
 \label{AovL}
\end{equation}
\label{thm:ocha}
\end{thm}
For an $A_\infty$-algebra $A$ over an $L_\infty$-algebra $L$, 
the substructure $(A,\{n_{0,k}\}_{k\ge 1})$ forms an $A_\infty$-algebra 
and the substructure $(L\oplus A,\{n_{p,1}\}_{p\ge 0})$ makes $A$ 
an $L_\infty$-module over $(L,\l)$. 
More precisely,
$n_{0,k}=m_k,\ n_{p,1} = k_{p+1}$; 
the $n_{1,q>0}$ map $L$ into ${\s}\Coder(T^c A)$
and the rest extend that to an $L_\infty$-map. 
This is just the higher homotopy structure 
a mathematician would construct 
by the usual procedures of strong homotopy algebra 
(see the Appendix by M.~Markl).

Here the sign exponent $\mu_{p,i}(\sigma)$ is given explicitly by 
\begin{equation}
\mu_{p,i}(\sigma)=\epsilon(\sigma)+
(c_{\sigma(1)}+\cdots +c_{\sigma(p)})+(o_1+\cdots+o_i)+
(o_1+\cdots +o_i)(c_{\sigma(p+1)}+\cdots+c_{\sigma(n)})\ ,
 \label{eta_pi}
\end{equation}
corresponding to the signs effected by the interchanges. 
The sign can be seen easily in the coalgebra and tree expressions. 
We can also write the defining equation (\ref{AovL}) 
in the following shorthand expression, 
\begin{equation*}
 \begin{split}
0=& \sum_{\sigma\in\S_{p+r=n}}(-1)^{\epsilon(\sigma)}
n_{1+r,m}\left( 
(l_p\otimes \1_c^{\otimes r}\otimes\1_o^{\otimes m})
(c_{\sigma(I)};o_1,\cdots,o_m) \right) \\
& +
\sum_{\substack{\sigma\in\S_{p+r=n}\\i+s+j=m}}
(-1)^{\epsilon(\sigma)}
n_{p,i+1+j}\left( (\1_c^{\otimes p}\otimes \1_o^{\otimes i}\otimes
n_{r,s}\otimes \1_o^{\otimes j} )
(c_{\sigma(I)};o_1,\cdots,o_m) \right)\ ,
 \end{split}
\end{equation*}
where the complicated sign is absorbed into this expression. 
Note that the rule for the action of tensor products of 
graded multi-linear maps on $(\cH_c)^{\otimes n}\otimes (\cH_o)^{\otimes m}$ 
is determined in a canonical way; for instance for 
$(f\otimes\cdots)(c_1,\cdots,c_n;o_1,\cdots,o_m)$ 
with the first multi-linear map 
$f:(\cH_c)^{\otimes k}\otimes (\cH_o)^{\otimes l}\raw\cH$, 
we may bring $(c_1,\cdots,c_k;o_1,\cdots,o_l)$ to the first $f$ 
with the associated sign, do the same thing for the next multi-linear 
map in $\cdots$ and repeat this in order.

String field theory suggests that an open-closed homotopy algebra 
includes the addition of the maps $n_{p,0}:L^{\otimes p}\raw A$ 
and in particular $n_{1,0}:L\raw A$ corresponding to 
the opening of a closed string to an open one. 
\begin{defn}[Open-Closed Homotopy Algebra (OCHA)]
An {\em open-closed homotopy algebra (OCHA)}
\footnote {The authors worked with the acronym for several weeks 
before realizing it is Japanese for `tea'.}
$(\cH=\cH_c\oplus \cH_o, \l, \n)$ consists of 
an $L_\infty$-algebra $(\cH_c,\l)$ and a family of maps 
$\n=\{n_{p,q}:\cH_c^{\ott p}\ott\cH_o^{\ott q} \to\cH_o\}$ 
for $p,q\geq 0$ with the exception of $(p,q)=(0,0)$ 
satisfying the compatability conditions (\ref{AovL}):
\begin{equation}
\begin{split}
&0=\sum_{\sigma\in\S_{p+r=n}}(-1)^{\epsilon(\sigma)}
n_{1+r,m}(l_p(c_{\sigma(1)},\cdots,c_{\sigma(p)}),
c_{\sigma(p+1)}\cdots,c_{\sigma(n)};o_1,\cdots,o_m) \\
& +
\sum_{\substack{\sigma\in\S_{p+r=n}\\i+s+j=m}}
(-1)^{\mu_{p,i}(\sigma)}
n_{p,i+1+j}(c_{\sigma(1)},\cdot\cdot,c_{\sigma(p)};o_1,\cdot\cdot,o_i,
n_{r,s}(c_{\sigma(p+1)},\cdot\cdot,c_{\sigma(n)};o_{i+1},\cdot\cdot,o_{i+s}),
o_{i+s+1},\cdot\cdot,o_m) \ , 
\end{split}
\end{equation}
for the full range $n,m\ge 0$, $(n,m)\ne (0,0)$. 
 
A {\em weak OCHA} consists of 
a weak $L_\infty$-algebra $(\cH_c,\l)$ with a family of maps 
$\n=\{n_{p,q}: \cH_c^{\ott p}\ott\cH_o^{\ott q} \to\cH_o\}$ 
now for $p,q\geq 0$ satisfying the analog of the above relation.
 \label{defn:ocha}
\end{defn}
For an OCHA $(\cH,\l,\n)$, 
the multi-linear maps $\{n_{p,q}\}_{p\ge 1,q\ge 0}$ still 
correspond to an 
adjoint 
$L_\infty$-map $\cH_c\raw\Coder(T^c\cH_o)$, 
 as in the case of an $A_\infty$-algebra 
over an $L_\infty$-algebra.
This has a particular importance 
in terms of deformation theory, cf. subsection \ref{ssec:defs}, 
where the addition of maps $n_{p,0}$ 
leads in turn to the deformation of 
the $A_\infty$-structure $\m$ to a {\em weak} $A_\infty$-structure.
\begin{defn}[Open-closed homotopy algebra (OCHA) morphism]
For two weak OCHAs $(\cH,\l,\n)$ and $(\cH',\l',\n')$, 
consider a collection $\f$ 
of degree zero (degree preserving) multi-linear maps 
\begin{equation*}
 \begin{split}
 f_k &: (\cH_c)^{\otimes k}
 \raw\cH'_c\ ,\qquad \mbox{for}\ {k\ge 0}\ , \\
 f_{k,l} &: (\cH_c)^{\otimes k}\otimes (\cH_o)^{\otimes l}
 \raw\cH'_o\ ,\qquad \mbox{for}\ {k,l\ge 0}\ ,
 \end{split}
\end{equation*}
where $f_k$ and $f_{k,l}$ are graded symmetric 
with respect to $(\cH_c)^{\otimes k}$. 
We call $\f:(\cH,\l,\n)\raw (\cH',\l',\n')$ 
a {\em weak OCHA-morphism} when 
$\{f_k\}_{k\ge 0}:(\cH_c,\l)\raw (\cH_c',\l')$ is 
a weak $L_\infty$-morphism and $\{f_{k,l}\}_{k,l\ge 0}$ 
further satisfies the following relations: 
\begin{equation}
\begin{split}
& \sum_{\sigma\in\S_{p+r=n}}(-1)^{\epsilon(\sigma)}
f_{1+r,m}\left( 
(l_p\otimes \1_c^{\otimes r}\otimes\1_o^{\otimes m})
(c_{\sigma(I)};o_1,\cdots,o_m) \right) \\
& +
\sum_{\substack{\sigma\in\S_{p+r=n}\\i+s+j=m}}
(-1)^{\epsilon(\sigma)}
f_{p,i+1+j}\left( (\1_c^{\otimes p}\otimes \1_o^{\otimes i}\otimes
n_{r,s}\otimes \1_o^{\otimes j} )
(c_{\sigma(I)};o_1,\cdots,o_m) \right) \\
&=
\sum_{\substack{\sigma\in\S_{(r_1+\cdots +r_i)+(p_1+\cdots +p_j)=n}\\
(q_1+\cdots +q_j)=m}}
\frac{(-1)^{\epsilon(\sigma)}}{i!}
n'_{i,j}\left( (f_{r_1}\otimes\cdot\cdot\otimes f_{r_i}\otimes
f_{p_1, q_1}\otimes\cdot\cdot\otimes f_{p_j,q_j}
)(c_{\sigma(I)};o_1,\cdots,o_m) \right) \ . 
\end{split}
 \label{ocmorpcd}
\end{equation}
The right hand side is written explicitly as 
\begin{equation*}
 \begin{split}
& n'_{i,j}\left( (f_{r_1}\otimes\cdots\otimes 
f_{r_i}\otimes
f_{p_1, q_1}\otimes\cdots\otimes 
f_{p_j,q_j})(c_{\sigma(I)};o_1,\cdots,o_m) \right) \\
&=(-1)^{\tau_{\oraw{p},\oraw{q}}(\sigma)}
n'_{i,j}(f_{r_1}
(c_{\sigma(1)},\cdot\cdot,c_{\sigma(r_1)}),
\cdots\cdot,f_{r_i}
(c_{\sigma({\bar r}_{i-1}+1)},\cdot\cdot,c_{\sigma({\bar r}_i)}); \\
&\hspace*{1.0cm}
f_{p_1, q_1}
(c_{\sigma({\bar r}_i+1)},\cdot\cdot,c_{\sigma({\bar p}_1)}
;o_1,\cdot\cdot,o_{q_1}),
\cdots\cdot,f_{p_j,q_j}
(c_{\sigma({\bar p}_{j-1}+1)},\cdot\cdot,c_{\sigma({\bar p}_j)};
o_{{\bar q}_{j-1}+1},\cdot\cdot,o_{{\bar q}_j})
)\ , 
 \end{split}
\end{equation*}
where 
${\bar r}_k:=r_1+\cdots +r_k$, 
${\bar p}_k:={\bar r}_i+ p_1+\cdots +p_k$, 
${\bar q}_k:=q_1+\cdots +q_k$ 
and $\tau_{\oraw{p},\oraw{q}}(\sigma)$ is given by 
$$
\tau_{\oraw{p},\oraw{q}}(\sigma)=
\sum_{k=1}^{j-1}
(c_{\sigma({\bar p}_k+1)}+\cdots +c_{\sigma({\bar p}_{k+1})})
(o_1+\cdots +o_{{\bar q}_k}) \ .
$$
In particular, if $(\cH,\l,\n)$ and $(\cH',\l',\n')$ are 
OCHAs and if $f_0=f_{0,0}=0$, 
we call it an {\em OCHA-morphism}. 
 \label{defn:ocmorp}
\end{defn}
\begin{defn}[OCHA-quasi-isomorphism]
Given two OCHAs 
$(\cH,\l,\n)$, 
$(\cH',\l',\n')$ and an OCHA-morphism 
$\f: (\cH,\l,\n)\raw (\cH',\l',\n')$, 
$\f$ is called an {\it OCHA-quasi-isomorphism} if 
$f_1+f_{0,1}:\cH\raw\cH'$ induces an isomorphism 
between the cohomology spaces of the complexes $(\cH,d:=l_1+n_{0,1})$ 
and $(\cH',d')$. 
In particular, if $f_1+f_{0,1}$ is an isomorphism, 
we call $\f$ an {\it OCHA-isomorphism}. 
 \label{defn:quasiisom}
\end{defn}

 \subsection{The coalgebra description}
\label{ssec:coalg}

Consider an OCHA $\cH=\cH_c\oplus\cH_o$. 
Recall that the separate $L_\infty$- and $A_\infty$-structures 
are described by coderivation differentials on, respectively, 
$C(\cH_c)$ and $T^c(\cH_o)$. 
The defining multi-linear maps for $\cH$ are to be extended 
to coderivations of $C(\cH_c)\ott T^c(\cH_o)$. 
The coproduct on $T^c(\cH_c)\ott T^c(\cH_o)$ is the standard tensor
product coproduct defined by 
\begin{equation}
 \begin{split}
  &\tri((c_1\otimes\cdots\otimes c_m)\otimes
 (o_1\otimes\cdots\otimes o_n)) \\
 &\quad =\sum_{p=0}^{m}
\sum_{q=0}^{n}
(-1)^{\eta(p,q)}(c_1\otimes\cdots\otimes c_p\otimes
o_1\otimes\cdots\otimes o_q)
\otimes (c_{p+1}\otimes\cdots\otimes c_m
\otimes o_{q+1}\otimes\cdots\otimes o_n)\ , 
 \end{split}
\end{equation}
\noindent
where $\eta(p,q)=(c_{p+1}+\cdots +c_m)(o_1+\cdots +o_q)$. 

The relevant subcoalgebra is $C(\cH_c)\ott T^c(\cH_o)$. 
Now we define the total coderivation $\l +\n$ by lifting
\begin{equation}
 \sum_{k\ge 1}(l_k+m_k) +\sum_{p\ge 1,q\ge 0}n_{p,q}\ ,
\end{equation}
with $m_k=n_{0,k}$. 
Thus we have an OCHA iff $\l +\n$ is a codifferential:
\begin{equation}
(\l +\n)^2=0\ .
 \label{occodiff}
\end{equation}
If this is true with the addition of $l_0$ and $m_0$, we have a weak OCHA.

Also, given two OCHAs $(\cH,\l,\n)$ and $(\cH',\l',\n')$, 
an OCHA-morphism $\f:(\cH,\l,\n)\raw(\cH',\l',\n')$ 
can be lifted to the coalgebra homomorphism 
$\f: C(\cH_c)\otimes T^c(\cH_o)\raw C(\cH'_c)\otimes T^c(\cH'_o)$ 
and the condition for an OCHA-morphism is written as 
$\f\circ (\l+\n)=(\l'+\n')\circ\f$.

 \subsection{The tree description}
\label{ssec:octree}

We associated the $k$-corolla of planar rooted trees to 
the multi-linear map $m_k$ of an $A_\infty$-algebra,   
and the $k$-corolla of non-planar rooted trees to the graded symmetric 
multi-linear map $l_k$ of an $L_\infty$-algebra. 
For an OCHA $(\cH,\l,\n)$, the corolla 
corresponding to $n_{k,l}$ should be expressed 
as the following mixed corolla, 
\begin{equation}
n_{k,l}\ \ \lglraw\ 
\begin{minipage}[c]{45mm}{\includegraphics{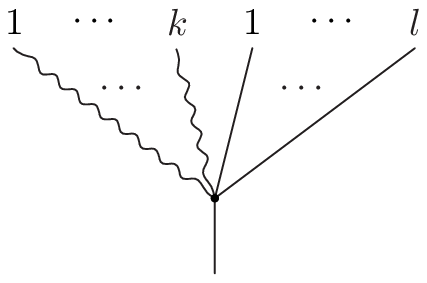}}
\end{minipage}\ ,
 \label{n_kl}
\end{equation}
which is partially symmetric (non-planar), that is, 
only symmetric with respect to the $k$ leaves. 
Let us consider such corollas for $2k+l+1\ge 3$ together with 
non-planar corollas $\{l_k\}_{k\ge 2}$. 
Since we have two kinds of edges, we have two kinds of grafting; 
grafting of edges associated to $\cH_c$ (closed string edges) and 
those for $\cH_o$ (open string edges). 
We denote them by $\circ_i$ and $\bullet_i$, respectively. 
For these corollas, we have three types of the composite; 
in addition to the composite $l_{1+k}\circ_i l_l$ in $\cL_\infty$, 
there is a composite $n_{k,m}\circ_i l_p$ described by 
\begin{equation}
\begin{minipage}[c]{45mm}{\includegraphics{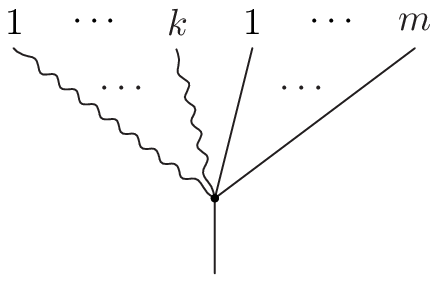}}
\end{minipage}
\circ_i
\begin{minipage}[c]{40mm}{\includegraphics{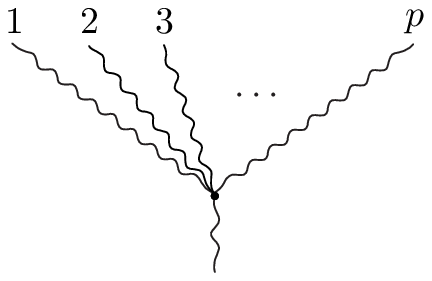}}
\end{minipage}\ =\ 
\begin{minipage}[c]{50mm}{\includegraphics{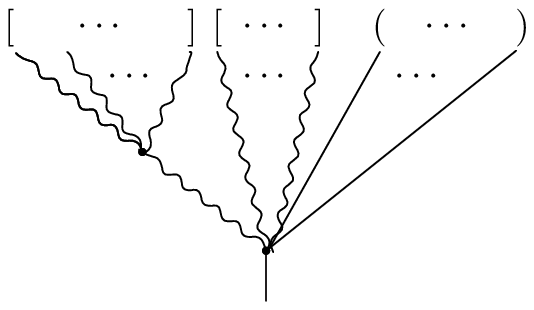}}
\end{minipage}\ ,
\end{equation}
where in the right hand side the labels are given by 
$[i,\cdots,i+p-1][1,\cdots,i-1,i+p,\cdots,p+k-1](1,\cdots,m)$, 
and the composite $n_{p,q}\bullet_i n_{r,s}$ 
\begin{equation}
\begin{minipage}[c]{45mm}{\includegraphics{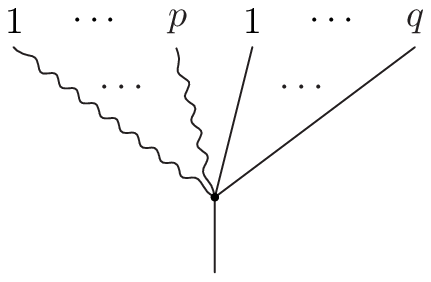}}
\end{minipage}
\bullet_i
\begin{minipage}[c]{40mm}{\includegraphics{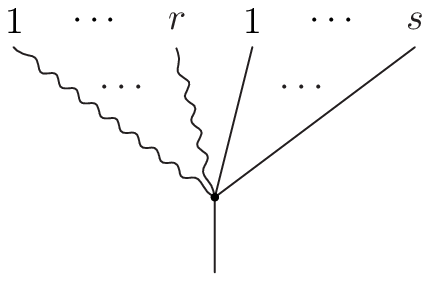}}
\end{minipage}\ =\ 
\begin{minipage}[c]{50mm}{\includegraphics{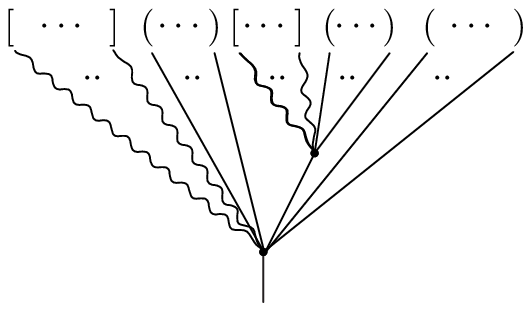}}
\end{minipage}\ 
\end{equation}
with labels 
$[1,\cdots,p](1,\cdots,i-1)[p+1,\cdots,p+r](i,\cdots,i+s-1)(i+s,\cdots,q+s-1)$.
To these resulting trees, grafting of a corolla $l_k$ or $n_{k,l}$ 
can be defined in a natural way, and we can repeat this procedure. 
Let us consider tree graphs obtained in this way, 
that is, by grafting the corollas $l_k$ and $n_{k,l}$ recursively, 
together with 
the action of permutations of the labels for closed string leaves. 
Each of them has a closed string root edge or 
an open string root edge. 
The tree graphs with closed string root edge, 
with the addition of the identity $e_c\in\cL_\infty(1)$, 
generate $\cL_\infty$ as stated in subsection \ref{ssec:Linftytree}. 
On the other hand, 
the tree graphs with open string root edge are new; 
the graded vector space generated by them with 
$k$ closed string leaves and $l$ open string leaves 
we denote by $\cN_\infty(k;l)$. 
In particular, we formally add the identity $e_o$ generating 
$\cN_\infty(0;1)$, and $\cN_\infty(1;0)$ is generated by a corolla
$n_{1,0}$. 
For $\cN_\infty:=\oplus_{k,l}\cN_\infty(k;l)$, 
the tree operad relevant here is then 
$\cOC_\infty:=\cL_\infty\oplus\cN_\infty$. 
For each tree $T\in\cOC_\infty$, 
its grading is given by the number of the vertices $v(T)$. 

{}For trees in $\cOC_\infty$, 
let $T'\raw T$ indicate that $T$ is obtained from $T'$ by 
contracting a closed or an open internal edge. 
A degree one differential $d: \cOC_\infty\raw\cOC_\infty$ 
is given by 
\begin{equation*}
  d(T)= \sum_{T'\raw T}\pm T'\ ,
\end{equation*}
so that the following compatibility holds: 
\begin{equation*}
 d(T\circ_i T')=d(T)\circ_i T' + (-1)^{v(T)} T\circ_i d(T')\ ,
 \qquad  d(T\bullet_i T'')=d(T)\bullet_i T'' 
 + (-1)^{v(T)} T\bullet_i d(T'')\ .
\end{equation*}
Thus, $\cOC_\infty$ forms a dg operad. 
In particular, $d(l_k)$ is given by eq.(\ref{d_Linfty-operad}), 
and $d(n_{n,m})$ is as follows: 
\begin{equation}
- \sum_{\sigma\in\S_{p+r=n}}
\begin{minipage}[c]{50mm}{\includegraphics{n_kll_p-num.eps}}
\end{minipage}
\ - 
\sum_{\substack{\sigma\in\S_{p+r=n}\\i+s+j=m}}
\begin{minipage}[c]{50mm}{\includegraphics{n_pq-rs-num.eps}}
\end{minipage}\ ,
 \label{occd-operad}
\end{equation}
where the labels for the first and the second terms are 
$[\ss(1),\cdots,\ss(p)] [\ss(p+1),\cdots,\ss(n)] (1,\cdots,m)$ and 
$[\ss(1),\cdots,\ss(p)] (1,\cdots,i) [\ss(p+1),\cdots,\ss(n)] 
(i+1,\cdots,i+s) (i+s+1,\cdots,m)$, respectively. 

An algebra $\cH:=\cH_c\oplus\cH_o$ over $\cOC_\infty$ is obtained by 
a representation 
\begin{equation*}
 \phi:\cL_\infty(k)\raw \Hom(\cH_c^{\otimes k},\cH_c) \ ,\qquad 
 \phi:\cN_\infty(k;l)\raw 
 \Hom((\cH_c)^{\otimes k}\otimes (\cH_o)^{\otimes l},\cH_o) \ 
\end{equation*}
which is compatible with respect to the grafting 
$\circ_i$, $\bullet_i$ and the differential $d$. 
Here, regarding elements in 
both $\Hom(\cH_c^{\otimes k},\cH_c)$ and 
$\Hom((\cH_c)^{\otimes k}\otimes (\cH_o)^{\otimes l},\cH_o)$ 
as those in $\Coder(C(\cH_c)\otimes T^c(\cH_o))$, 
the differential in the algebra side is given by 
$[l_1+n_{0,1},\ \, ]$. 
By combining it with eq.(\ref{occd-operad}) 
one can recover the condition of an OCHA (\ref{AovL}). 

If we adjust the notation for grading as 
${\downarrow}{\downarrow}\cH_c$ and ${\downarrow}\cH_o$, 
the degree of the multi-linear map $l_k$ is $3-2k$ as stated previously 
and the degree of $n_{k,l}$ turns out to be $1+(1-l)-2k=2-(2k+l)$. 
The grading of a tree $T\in\cN_\infty(k;l)$ is then replaced by 
$int(T)+(2-2k-l)$, which is equal to minus the dimension of the 
corresponding boundary piece of the compactified moduli space 
of a disk with $k$ points interior and $l$ points on the boundary 
(see \cite{Z2}).

\section{Cyclic structures}
\label{sec:cyclic}

Now we consider an additional structure, {\em cyclicity}, 
on open-closed homotopy algebras.  
Algebras with invariant inner products ($\la ab,c\ra = \la a,bc\ra$ or
$\la[a,b],c\ra = \la a,[b,c]\ra$) are very important in mathematical physics;
the analogous definition for strong homotopy algebras is straightforward
(cf. \cite{MSS}, sections II.5.1 and II.5.2). 
The string theory motivation for this additional structure is 
that punctures on the boundary
of the disk inherit a cyclic order from the orientation of the disk
and the operations are to respect this cyclic structure, just as the
$L_\infty$-structure reflects the symmetry of the punctures 
in the interior of the disk or on the sphere.

In our context, cyclicity is defined 
in terms of constant symplectic inner products.
(The terminology is that used for symplectic structures on
supermanifolds \cite{AKSZ}; 
see also \cite{thesis} and references therein. 
These inner products are 
also essential to the description of the Lagrangians 
appearing in string field theory.)

\begin{defn}[Constant symplectic structure]
Bilinear maps, 
$\omega_c: \cH_c\otimes\cH_c\raw\C$ and 
$\omega_o: \cH_o\otimes\cH_o\raw\C$, 
are called {\em constant symplectic structures} 
when they have fixed integer degrees $|\omega_c|, |\omega_o|\in\Z$ 
and are non-degenerate and skew-symmetric. 
Here `skew-symmetric' indicates that 
\begin{equation*}
 \omega_c(c_2,c_1)=-(-1)^{c_1c_2}\omega_c(c_1,c_2)\ ,\qquad 
 \omega_o(o_2,o_1)=-(-1)^{o_1o_2}\omega_o(o_1,o_2)\ 
\end{equation*} 
for any $c_1,c_2\in\cH_c$, $o_1,o_2\in\cH_o$, and 
degree $|\omega_c|$ and $|\omega_o|$ implies that 
$\omega_c(c_1,c_2)=0$ except for $\deg(c_1)+\deg(c_2)+|\omega_c|=0$ and
$\omega_o(o_1,o_2)=0$ except for $\deg(o_1)+\deg(o_2)+|\omega_o|=0$. 
We further denote the constant symplectic structure on 
$\cH=\cH_c\oplus\cH_o$ by $\omega:=\omega_c\oplus\omega_o$. 
 \label{defn:sym}
\end{defn}
Suppose that an OCHA $(\cH,\l,\n)$ 
is equipped with constant symplectic structures 
$\omega_c: \cH_c\otimes\cH_c\raw\C$ and 
$\omega_o: \cH_o\otimes\cH_o\raw\C$ as in Definition \ref{defn:sym}. 
For $\{l_k\}_{k\ge 1}$ and $\{n_{p,q}\}_{p+q\ge 1}$, 
let us define two kinds of 
multi-linear maps by 
\begin{equation*}
 \V_{k+1}=\omega_c(l_k\ott\1_c): (\cH_c)^{\otimes (k+1)}\raw\C\
 ,\qquad 
 \V_{p,q+1}=\omega_o(n_{p,q}\ott\1_o):(\cH_c)^{\otimes p}\otimes 
(\cH_o)^{\otimes (q+1)}\raw\C\ 
\end{equation*}
or more explicitly 
\begin{equation*}
 \V_{k+1}(c_1,\cdots,c_{k+1})
=\omega_c(l_k(c_1,\cdots,c_k),c_{k+1}) 
\end{equation*}
and 
\begin{equation*}
 \V_{p,q+1}(c_1,\cdots,c_p;o_1,\cdots,o_{q+1})=
 \omega_o(n_{p,q}(c_1,\cdots,c_p;o_1,\cdots,o_{q}),o_{q+1})\ .
\end{equation*}
The degree of $\V_{k+1}$ and $\V_{p,q+1}$ are 
$|\omega_c|+1$ and $|\omega_o|+1$. 
\begin{defn}[Cyclic open-closed homotopy algebra (COCHA)]
An OCHA $(\cH,\omega,\l,\n)$ 
is a {\em cyclic open-closed homotopy algebra} (COCHA)
when 
$\V_{k+1}$ is graded symmetric with respect to 
any permutation of $(\cH_c)^{\otimes(k+1)}$ and 
$\V_{p,q+1}$ has cyclic symmetry 
with respect to cyclic permutations of $(\cH_o)^{\otimes (q+1)}$, 
that is, if 
\begin{equation*}
 \V_{k+1}(c_1,\cdots,c_{k+1})
 =(-1)^{\epsilon(\sigma)}\V_{k+1}(c_{\sigma(1)},\cdots,c_{\sigma(k+1)})\ , 
 \qquad \sigma\in\S_{k+1}\ 
\end{equation*}
and 
\begin{equation*}
 \V_{p,q+1}(c_1,\cdots,c_p;o_1,\cdots,o_{q+1})=
(-1)^{o_1(o_2+\cdots o_{q+1})}
\V_{p,q+1}(c_1,\cdots,c_p;o_2,\cdots,o_{q+1},o_1)\ . 
\end{equation*} 
The graded commutativity of $\V_{p,q+1}$ with respect to permutations of
$(\cH_c)^{\ott p}$, that is, 
\begin{equation*}
 \V_{p,q+1}(c_1,\cdots,c_p;o_1,\cdots,o_{q+1})=
 (-1)^{\epsilon(\sigma)}\V_{p,q+1}
(c_{\sigma(1)},\cdots,c_{\sigma(p)};o_1,\cdots,o_{q+1})\ ,
 \qquad \sigma\in\S_p
\end{equation*}
automatically holds by the definition of $\n$. 
 \label{defn:cyc} 
\end{defn}
Since we have non-degenerate inner products $\omega_c$ and $\omega_o$,
we can identify $\cH$ with its linear dual, 
then reverse the process and define further maps
$$
r_{p-1,q+1}: (\cH_c)^{\ott (p-1)}\ott (\cH_o)^{\ott (q+1)}\raw \cH_c
$$
with relations amongst themselves and with the operations already 
defined, which can easily be deduced from their definition.
In particular, for $n_{1,0}:\cH_c\raw\cH_o$ we have 
$r_{0,1}: \cH_o\raw\cH_c$.  
Namely, for the cyclic case 
the fundamental object is the multi-linear map $\V_{p,q+1}$ 
where $n_{p,q}$ and $r_{p-1,q+1}$ are equivalent under the 
relation above. 
However, we get a codifferential (\ref{occodiff}) 
since we took $n_{p,q}$ instead of $r_{p-1,q+1}$ for defining an OCHA. 
Physically, for the multi-linear map $\V_{p,q+1}$, 
choosing $\cH_o$ as a root edge instead of $\cH_c$ 
as in eq.(\ref{n_kl}) 
is related to a standard compactification of the corresponding 
Riemann surface (a disk with $p$ points interior and $(q+1)$ points 
on the boundary).

\section{Minimal model theorem and decomposition theorem}
\label{sec:MMth}

Homotopy algebras are designed to  have  homotopy invariant properties. 
A key and useful theorem in homotopy algebras is then the minimal model
theorem. For $A_\infty$-algebras, it was proved by 
Kadeishvili \cite{kadei1}. 
For the construction of minimal models of $A_\infty$-structures, 
in particular on the homology of a differential graded algebra, 
homological perturbation theory (HPT) is developed by 
\cite{gugen,hueb-kadei,GS,gugen-lambe,GLS:chen,GLS}, 
for instance, and the form of a minimal model is also given explicitly 
and more recently in \cite{mer,KS}. 

There are various results referred to as minimal model theorems:
the weakest form asserts the existence of a quasi-isomorphism as 
$A_\infty$-algebras 
$H(A) \to A$ for an $A_\infty$-structure on $H(A)$,
by noticing that all the relevant obstructions
vanish because the homology of $A$ and $H(A)$ agree.
A stronger result constructs an $A_\infty$-structure on $H(A)$ and
the quasi-isomorphism,
then  a decomposition  theorem is proved from which
the inverse quasi-isomorphism follows \cite{thesis,Le-Ha,KaTe,kadei1}.
Alternatively, the full strength of homological perturbation theory 
gives the maps in
both directions and the homotopy for the composition $A\to H(A) \to A$
all together.

The corresponding theorems for $L_\infty$-algebras are more recent: 
\cite{schuhmacher} for the two step procedure, 
\cite{jh-jds} for the full HPT treatment. 
The latter points out that, although $L_\infty$-algebras can be 
constructed by symmetrization of $A_\infty$-algebras, 
the corresponding constructions of the maps and homotopy
are more subtle.

It is not surprising that the minimal models and decompositions 
exist also for our OCHAs. 
These theorems imply that, for an OCHA $(\cH,\l,\n)$, 
the higher multi-linear structures 
$l_k$, $k\ge 2$ and $n_{p,q}$, $(p,q)\ne (0,1)$ have been transformed
to those on $H(\cH)$, 
where $H(\cH)$ is the cohomology of the complex $(\cH,d=l_1+n_{0,1})$. 
Even though some of those higher structures 
may have been zero on the original OCHA $\cH$, 
those on $H(\cH)$ need not be.

We present these statements more precisely below, leaving detailed
proofs to the industrious reader. 
\begin{defn}[Minimal open-closed homotopy algebra]
An OCHA $(\cH=\cH_c\oplus\cH_o,\l,\n)$ 
is called {\it minimal} if 
$l_1=0$ on $\cH_c$ and $n_{0,1}=0$ on $\cH_o$. 
 \label{defn:minimalalg}
\end{defn}
\begin{defn}[Linear contractible open-closed homotopy algebra]
A {\em linear contractible} OCHA $(\cH,\l,\n)$ is a 
complex $(\cH,d=l_1+n_{0,1})$ which has trivial cohomology, 
that is, an OCHA $(\cH=\cH_c\oplus\cH_o,\l,\n)$ 
such that $l_l=0$ for $l\ge 2$, 
$n_{p,q}=0$ except for $(p,q)=(0,1)$, and the complexes 
$(\cH_c,l_1)$, $(\cH_o,n_{0,1})$ having trivial cohomologies. 
 \label{defn:cont}
\end{defn}
\begin{thm}[Decomposition theorem for open-closed homotopy algebras]
Any OCHA is 
isomorphic to the direct sum of a minimal OCHA 
and a linear contractible OCHA. 
 \label{thm:MandC}
\end{thm}
A weak version of the minimal model theorem follows from the
decomposition theorem above: 
\begin{thm}[Minimal model theorem for open-closed homotopy algebras]
{}For a given OCHA $(\cH,\l,\n)$, 
there exists a minimal OCHA $(H(\cH),\l',\n')$ and 
an OCHA-quasi-isomorphism $\f: (H(\cH),\l',\n')\raw (\cH,\l,\n)$. 
In particular, the minimal model can be taken so that 
$l'_2=H(l_2)$, $n'_{0,2}=H(n_{0,2})$ and $n'_{1,0}=H(n_{1,0})$. 
 \label{thm:minimal}
\end{thm}
To obtain a homotopy equivalence 
from an initial quasi-isomorphism $\f$ above, 
one way is to employ the decomposition theorem 
(Theorem \ref{thm:MandC}). 
Alternatively, it 
can be obtained directly by the methods of HPT 
(see Theorem \ref{thm:minimalHPT} below).
In either approach, 
for a given OCHA $(\cH,\l,\n)$, 
one first considers a Hodge decomposition of the 
complex $(\cH,d=l_1+n_{0,1})$. 
Namely, decompose $\cH$ into a direct sum isomorphic to
$\cH=H(\cH)\oplus\cC$, $\cC:=Y\oplus dY$ 
with a contracting homotopy $h:dY\raw Y$ of degree minus one. 
Together with the inclusion $\iota$ and the projection $\pi$, 
let us express these data as 
\begin{equation*}
 \Nsddata {H(\cH)} {\hspace*{0.3cm}\iota}{\hspace*{0.3cm}\pi}{\cH}h \ .
\end{equation*}

Then the decomposition theorem (Theorem \ref{thm:MandC}) states that 
there exists an OCHA-isomorphism 
$\f_{isom}: (H(\cH),\l',\n')\oplus (\cC,d)\raw (\cH,\l,\n)$, where 
$(\cC,d)$ is the linear contractible OCHA. 
The OCHA-isomorphism is obtained by first decomposing 
the $L_\infty$-algebra $(\cH_c,\l)$ 
into the direct sum of a minimal part $H(\cH_c)$ 
and a linear contractible part $\cC_c$, and then decomposing 
the OCHA $(\cC_c,d_c)\oplus (H(\cH_c)\oplus\cH_o,\l',\n)$ in a 
similar way as in the $A_\infty$ case. 
Because we have the OCHA-isomorphism, 
we may consider a homotopy equivalence 
between $(H(\cH),\l',\n')$ and $(H(\cH),\l',\n')\oplus (\cC,d)$. 
In fact, the maps $\iota$ and $\pi$ naturally extend to 
OCHA-quasi-isomorphisms between them, and 
the corresponding homotopy is obtained as in the sense in 
Theorem \ref{thm:minimalHPT} below 
(see \cite{thesis} for the $A_\infty$ case) 
or as a path between them with some appropriate compatibility 
(\cite{KaTe} for the $A_\infty$ case). 

Alternatively, one can refine the 
standard HPT machinery to function in the category of OCHAs and their
morphisms or apply the known results to the $L_\infty$-algebra $\cH_c$ 
and then extend to $\cH$, 
regarding $\cH_o$ as an analog of an sh-algebra over $\cH_c$. 
The extra detail of the HPT form of the minimal model theorem is then:
\begin{thm} [HPT minimal model theorem for open-closed homotopy algebras]
\hfill\\
Given an OCHA $(\cH,\l,\n)$ and a Hodge decomposition with a contraction
\begin{equation*}
 \Nsddata {H(\cH)} {\hspace*{0.3cm}\iota}{\hspace*{0.3cm}\pi}{\cH}h \ ,
\end{equation*}
the linear maps $\pi$ and $\iota$ can be extended to coalgebra maps
and perturbed so that there exists a corresponding contraction of coalgebras
\begin{equation}
\Nsddata {\, \ C(H(\cH_c))\otimes
T^c(H(\cH_o))}{\hspace*{0.66cm}{\bar\iota}}{\hspace*{0.66cm}{\bar\pi}}
{C(\cH_c)\otimes T^c(\cH_o)}{\ \bar h\ \ }\ , 
 \label{HPT}
\end{equation}
where $\bar h$ is a degree minus one linear homotopy on 
$C(\cH_c)\otimes T^c(\cH_o)$, not necessarily a coalgebra map.
\label{thm:minimalHPT}
\end{thm}
In the same way as in the case of $A_\infty$-algebras, 
the minimal model theorem 
together with these additional theorems 
implies various corollaries. For instance, 
\begin{cor}[Uniqueness of minimal open-closed homotopy algebras]
For an OCHA $(\cH,\l,\n)$, 
its minimal OCHA $H(\cH)$ is unique up to an isomorphism on $H(\cH)$. 
\label{cor:uniqueness}
\end{cor}
\begin{cor}[Existence of an inverse quasi-isomorphism]
For two OCHAs
$(\cH,\l,\n)$ and $(\cH',\l',\n')$, 
suppose there exists an OCHA quasi-isomorphism 
$\f:(\cH,\l,\n)\raw(\cH',\l',\n')$. 
Then, there exists an inverse 
OCHA quasi-isomorphism 
$\f^{-1}:(\cH',\l',\n')\raw(\cH,\l,\n)$.  
 \label{cor:inverse}
\end{cor}
In particular, Corollary \ref{cor:inverse} 
guarantees that quasi-isomorphisms do in fact define a (homotopy) 
equivalence relation and in addition give bijective maps between 
the moduli spaces of the solution space 
of the corresponding Maurer-Cartan equations for 
quasi-isomorphic sh-algebras (see Theorem \ref{thm:moduli}). 

The same facts should hold also for cyclic OCHAs.

\section{Deformations and moduli spaces of $A_\infty$-structures}
\label{sec:deform}

 \subsection{Deformations and Maurer-Cartan equations}
\label{ssec:defs}

Consider
an OCHA $(\cH=\cH_c\oplus\cH_o,\l,\n)$. 
We will show how the combined structure implies 
the $L_\infty$-algebra $(\cH_c,\l)$ controls some
deformations of the $A_\infty$-algebra
$(\cH_o,\{m_k\}_{k\ge 1})$. We will further investigate the deformations
of this control as $\cH$ is deformed.

We first review some of the basics of deformation theory from a 
homotopy algebra point of view. 
The philosophy  of deformation theory which we follow 
(due originally, we believe, to Grothendieck
\footnote{See \cite{doran:bib} 
for an extensive annotated bibliography of deformation theory.}
 cf. \cite{SS, goldman-millson,deligne}) 
regards any deformation theory as `controlled' 
by a dg Lie algebra $\g$ 
(unique up to homotopy type as an $L_\infty$-algebra).

For the deformation theory of an (ungraded) associative algebra $A$, 
the standard controlling dg Lie algebra is $\Coder(T^c A)$ 
with the graded commutator 
as the graded Lie bracket \cite{jds:intrinsic}.  
Under the identification (including a shift in grading)
of $\Coder (T^c A)$ with $\Hom(T^cA, A)$
(which is the Hochschild cochain complex), 
this bracket is identified with the Gerstenhaber bracket 
and the differential with the Hochschild differential,
which can be written as $[m, \ ]$ \cite{gerst:coh}. 

The generalization to a differential graded associative algebra is 
straightforward; the differential is now: $[d_A+m_2,\quad]$. 
For an $A_\infty$-algebra,
the differential similarly generalizes to $[\m,\ ]$. 

Deformations of $A$ correspond to certain elements of  
$\Coder (T^c A)$, namely 
those that are solutions of an {\em integrability} equation, 
now known more commonly as a {\em Maurer-Cartan} equation.
\begin{defn}[The classical Maurer-Cartan equation] 
In a dg Lie algebra $(\g,d,[\ ,\ ])$,
the {\em classical  Maurer-Cartan equation} is
\begin{equation}
 d\theta + \ov{2}[\theta,\theta]=0
 \label{mceq-dgla}
\end{equation}
for $\theta\in \g^1 =( \downarrow L)^1$. 
\end{defn}
For an $A_\infty$-algebra $(A,\m)$ and $\theta \in \Coder^1(T^c A)$, 
a deformed $A_\infty$-structure 
is given by $\m + \theta$ iff 
\begin{equation*}
 (\m + \theta)^2 = 0\ .
\end{equation*}
Teasing this apart, since we start with $\m^2=0$, we have equivalently
\begin{equation}
 D\theta + 1/2[\theta, \theta] = 0\ ,
 \label{mceqCoder}
\end{equation}
hence the Maurer-Cartan name. (Here $D$ is the natural differential on
$\Coder(T^c A)\subset End (T^c A)$, i.e. $D\theta = [\m,\theta]$. )
Notice that we call this the Maurer-Cartan equation 
for the dg Lie algebra $(\Coder(T^cA), D, [\ ,\ ])$ 
but not for $(A,\m)$. 

For $L_\infty$-algebras, the analogous remarks hold, substituting
the Chevalley-Eilenberg complex for that of Hochschild, i.e. using
$\Coder\ C(L) \simeq \Hom (C(L), L)$. 

In any case, 
formal deformation theory controlled by a dg Lie algebra 
$(\g,d,[\ ,\ ])$ proceeds as follows. 
Consider a formal solution $\theta$ of the Maurer-Cartan equation 
(\ref{mceq-dgla}) in $\theta\in \g^1\otimes \hbar\C[[\hbar]]$, 
where $\hbar$ is a formal parameter. 
We express it as $\theta=\theta_{(1)}\hbar+\theta_{(2)}\hbar^2+\cdots$, 
where $\theta_{(i)}\in\g^1$. The Maurer-Cartan equation holds separately 
in different powers of $\hbar$, so we have 
\begin{align}
 (\hbar)^1 : & \qquad d\theta_{(1)}=0\ , \label{fdeform1}\\
 (\hbar)^2 : & \qquad d\theta_{(2)}
 +\half[\theta_{(1)},\theta_{(1)}]=0\ ,\label{fdeform2} \\
 (\hbar)^3 : & \qquad d\theta_{(3)}
 +[\theta_{(1)},\theta_{(2)}] =0\ , \label{fdeform3}\\
  \cdots \ \ & \qquad\qquad \cdots\cdots\qquad\qquad .  \nn
\end{align}
The first order solution $\theta_{(1)}$ is defined by the first equation 
(\ref{fdeform1}), that is, $\theta_{(1)}$ is a cocycle. 
This is also known as an infinitesimal deformation. 
We may proceed to second order 
if there is some $\theta_{(2)}$ satisfying the second equation 
(\ref{fdeform2}). 
Similarly, we can ask for $\theta_{(3)}$ satisfying the 
third equation (\ref{fdeform3}), etc. .

Since deformation theory is controlled by a dg Lie algebra 
{\em up to homotopy} (see Theorem \ref{thm:known}), 
the Maurer-Cartan equation should be extended to 
that for an $L_\infty$-algebra. 
We present the definition in the suspended $(L={\s}\g)$ notation. 

In addition to the convergence problem which would occur 
in the dg Lie algebra case, 
for an $L_\infty$-algebra on $L$ 
the Maurer-Cartan equation itself 
does not make sense in general 
since it consists of an infinite sum (see below). 
One way to avoid these problems is again to consider formal
deformation theory; 
one usually considers a homotopy algebra on a graded vector space $V$ 
over $\hbar\C[[\hbar]]$, or more generally 
a finite dimensional nilpotent commutative associative algebra. 
In particular, for an Artin algebra $A$ and its maximal ideal $m_A$, 
the standard way is to consider $V\otimes m_A$, where 
the degree of $A$ is set to be zero. 
The multi-linear operations on $V$ are extended 
to those on $V\otimes m_A$ trivially. 
{}From now on, we shall assume but not mention explicitly 
that any homotopy algebra $V$ we consider has been tensored with $m_A$ 
for some fixed $m_A$ and denote the result also by $V$. 
\begin{defn}
[The strong homotopy Maurer-Cartan equation] 
In an $L_\infty$-algebra $(L,\l)$, 
the {\em (generalized) Maurer-Cartan equation} is
\begin{equation*}
 \sum_{k\geq 1} \ov{k!} l_k(\cb,\cdots,\cb) = 0
\end{equation*}
for $\cb\in L^0$. 
\end{defn}
Note that the degree of $\cb$ is zero since $\g^1=L^0$. 
We denote the set of solutions of the Maurer-Cartan equation as 
$\cMC(L,\l)$. 
In the same sense as in the dg Lie algebra case, 
a cocycle $\cb\in L^0$, $l_1\cb=0$ play the role of 
a first order solution. 

Recall that for an OCHA $(\cH,\l,\n)$,
the adjoints of the maps $n_{p,q}$ constitute an $L_\infty$-map
\begin{equation}
\rho:\cH_c\raw {\s}\Coder(T^c\cH_o)\ .
\end{equation}
Since it is known that an $L_\infty$-morphism preserves 
the solutions of the Maurer-Cartan equations, we obtain the following: 
\begin{thm}
If $\cb\in\cH_c$ is a Maurer-Cartan element, 
then $\rho(\cb)\in {\s}\Coder(T^c\cH_o)$ gives a deformation 
of $\cH_o$ as a {\em weak} $A_\infty$-algebra.
 \label{thm:defs}
\end{thm}
In particular, a first order solution for $\cb\in\cH_c$ 
is preserved to be a 
first order solution in ${\s}\Coder(T^c\cH_o)$ by an $L_\infty$-morphism. 
The corresponding situation is the chain map (\ref{rho-chain}) 
considered previously: 
\begin{equation*}
 \rho(d_\g({\s}X)) = [\m, \rho({\s}X)]\ ,
\end{equation*}
where ${\s}X\in\cH_c$ and $\rho({\s}X)\in{\s}\Coder(T^c\cH_o)$. 
If ${\s}X$ is a first order solution, 
the chain map gives us $[\m,\rho({\s}X)]=0$. 
However, notice that this $\rho({\s}X)$, a first order solution in 
${\s}\Coder(T^c\cH_o)$, 
in general includes a constant term $\C\raw\cH_o$ coming from $n_{1,0}$. 
Namely, the first order deformation of $\m$ turns out to be 
a `{\em weak}' homotopy derivation, 
a natural extension of a strong homotopy derivation 
in Definition \ref{defn:shder} 
by including a map $\theta_0:\C\raw A$. 

Rather than treat this result in isolation, 
we look at more general deformations of $\cH$ as an OCHA. 
In order to do it, let us first explain another aspect of the Maurer-Cartan 
equation for a dg Lie algebra or an $L_\infty$-algebra more generally. 
\begin{lem}
For an $L_\infty$-algebra $(L,\l)$ and a graded 
vector space $L'$, 
consider a coalgebra isomorphism $\f: C(L')\raw C(L)$, that is, 
a collection of degree zero graded symmetric maps $\{f_0,f_1,\cdots \}$ 
such that $f_1:L'\raw L$ is an isomorphism. 
Then, the inverse of $\f$ exists, and 
a unique weak $L_\infty$-structure $\l'$ is induced by 
$\l'=(\f)^{-1}\circ\l\circ\f$ so that 
$\f:(L',\l')\raw (L,\l)$ is a weak $L_\infty$-isomorphism. 
 \label{label:weakL}
\end{lem}
It is clear by definition that $\l'$ is a degree one coderivation 
and $(\l')^2=0$. 

Moreover, if we take $\{f_0=\cb\in L^0, f_1=\1, f_2=\cdots=0 \}$ 
for $\f$, the explicit form of $\l'$ is given as follows 
(see Getzler \cite{getzler-L} and Schuhmacher \cite{schuhmacher}): 
\begin{equation}
 l'_l(c_1,\cdots,c_l):=\sum_{n\ge 0}\ov{n!}
 l_{n+l}(\cb^{\otimes n},c_1,\cdots,c_l)\ ,\qquad l\ge 0\ .
 \label{weakL}
\end{equation}
Here recall that $f_0:\C \to \cH_c$ so we identify $f_0$
with its image $\bar c$. 
Notice that $l'_0=\sum_{k\ge 1}\ov{k!}l_k(\cb^{\otimes k})$. 
Thus, $\l'$ gives a (strict) $L_\infty$-structure iff 
$\cb\in\cMC(L,\l)$. 

In this argument, we can also begin with a weak $L_\infty$-algebra 
$(L,\l)$ together with a straightforward modification of the 
Maurer-Cartan equation for a weak $L_\infty$-algebra. 

The same fact holds true also for (weak) $A_\infty$-algebras, 
as explained in subsection 2.4 in \cite{thesis} 
(the explicit form of the deformed $A_\infty$-structures 
can be found in \cite{Fukaya2,FOOO}). 

Let us consider the same story for an OCHA. 
\begin{lem}
For an OCHA $(\cH,\l,\n)$ and a graded vector space $\cH'$, 
consider a coalgebra map 
$\f: C(\cH'_c)\otimes T^c(\cH'_o)\raw C(\cH_c)\otimes T^c(\cH_o)$ 
such that $f_1+f_{0,1}:\cH'\raw\cH$ is an isomorphism. 
Then, a unique weak OCHA-structure $\l'+\n'$ is induced by 
$\l'+\n'=(\f)^{-1}\circ(\l+\n)\circ\f$ so that 
$\f:(\cH',\l',\n')\raw (\cH,\l,\n)$ is a weak OCHA-isomorphism. 
 \label{lem:weakOCHA}
\end{lem}
Again, the fact that $\f$ is a coalgebra map and $\l+\n$ is a 
degree one coderivation implies that $\l'+\n'$ is in fact a 
degree one coderivation, and $(\l'+\n')^2=0$ follows from 
$(\l+\n)^2=0$. 
The reason the structure is {\em weak} is the presence 
of the operations $l'_0$ and $n'_{0,0}$.

In particular, when we take a weak OCHA-isomorphism $\f$ given by 
\begin{equation*}
 f_0=\cb\in\cH_c^0\ ,\quad f_{0,0}=\ob\in\cH_o^0\ ,\quad 
f_1=\1_c\ ,\quad f_{0,1}=\1_o
\end{equation*}
and other higher multi-linear maps set to be zero, 
the deformed weak OCHA structure is given by $\l'$ 
in eq.(\ref{weakL}) and 
\begin{equation}
 \begin{split}
& n'_{p,q}(c_1,\cdots,c_p;o_1,\cdots,o_q)\\
 &\ \ :=\sum_{n,m_0,\cdots,m_k\ge 0}
 \ov{n!}n_{n+p,m_0+\cdots +m_k+q}
 (\cb^{\otimes n},c_1,\cdots,c_p;
 \ob^{\otimes m_0},o_1,\ob^{\otimes m_1},\cdots,
 \ob^{\otimes m_{q-1}},o_q,\ob^{\otimes m_q}) 
 \end{split}
 \label{ocha'}
\end{equation}
for $p\ge 0$ and $q\ge 0$. 

Now, we can spell out 
generalized Maurer-Cartan equations for OCHAs. 
\begin{defn}[Maurer-Cartan equations for $(\cH,\l,\n)$] 
For an OCHA $(\cH,\l,\n)$ 
and degree zero elements $\cb\in\cH_c$ and $\ob\in\cH_o$, 
we define
\begin{equation}
 \l_*(\cb):=\sum_k \ov{k!}l_k(\cb,\cdots,\cb)\ ,\qquad 
 \n_*(\cb;\ob):=\sum_{k,l}\ov{k!} 
 n_{k,l}(\cb,\cdots,\cb;\ob,\cdots,\ob)\ .
 \label{mcelt}
\end{equation}
We call the following pair of equations 
\begin{equation}
 0=\l_*(\cb)\ ,\qquad  0=\n_*(\cb;\ob)
 \label{mceq}
\end{equation}
the {\em Maurer-Cartan equations} 
for the OCHA $(\cH,\l,\n)$. 

The solution space of the Maurer-Cartan equations is denoted by 
\begin{equation*}
 \cMC(\cH,\l,\n)=\{(\cb,\ob)\in (\cH_c^0,\cH_o^0)\ |\ 
 \l_*(\cb)=0 ,\ \n_*(\cb;\ob)=0\ \}\ .
\end{equation*}
 \label{defn:mceq}
\end{defn}
The Maurer-Cartan equations (\ref{mceq}) 
are nothing but the condition that $l'_0=0$ and $n'_0=0$, 
since $l'_0=\l_*(\cb)$ and $n'_0=\n_*(\cb;\ob)$. 
In particular, the first equation is just 
the Maurer-Cartan equation for the $L_\infty$-algebra $(\cH_c,\l)$.

Now, one gets the following. 
\begin{thm}[Maurer-Cartan elements as deformations]
$(\cb,\ob)\in\cMC(\cH,\l,\n)$ gives a deformation of $(\cH_o,\m)$ 
as a (strict) $A_\infty$-algebra. 
 \label{thm:defAinfty}
\end{thm}
The explanations are as follows. 
First of all, for a weak OCHA $(\cH'=\cH_c'\oplus\cH'_o,\l',\n')$ 
given in eq.(\ref{weakL}) and eq.(\ref{ocha'}), 
let us consider its restriction to $\cH'_c=0$, that is, 
consider the defining equation for a (weak) OCHA (\ref{AovL}) 
and set $c_1=\cdots=c_n=0$. 
Then, only the equations for $n=0$ survive, which are given by 
\begin{equation*}
0=n'_{1,m}(l'_0;o_1,\cdots,o_m) 
+\sum_{i+s+j=m}
(-1)^{\beta(s,i)}
n'_{0,i+1+j}(\emptyset;o_1,\cdot\cdot,o_i,
n'_{0,s}(\emptyset;o_{i+1},\cdot\cdot,o_{i+s}),
o_{i+s+1},\cdot\cdot,o_m) \ . 
 \label{AovL'}
\end{equation*}
Here $l'=0$ iff $\cb\in\cMC(\cH_c,\l)$, then 
the first term in the right hand side drops out and 
the second term turns out to be the defining equation for 
a weak $A_\infty$-algebra. 
This is just the situation of Theorem \ref{thm:defs}, where 
$\rho(\cb)$ is given explicitly by 
\begin{equation*}
 {\downarrow}\rho(\cb)
 =\sum_{p\ge 1,q\ge 0}\ov{p!}n_{p,q}(\cb^{\otimes p};\ ,\cdots,\ )\ 
 \in \Hom(T^c\cH_o,\cH_o)\ .
\end{equation*} 
Note that, $(\cb, 0)\in(\cH^0_c,\cH^0_o)$ need not belong to 
$\cMC(\cH,\l,\n)$ even if $\cb\in\cMC(\cH,\l)$
because of the existence of $n_{k,0}(\cb,\cdots,\cb)$ terms 
in the second equation in eq.(\ref{mcelt}). 
Alternatively, for $\cb\in\cMC(\cH_c,\l)$, 
if we can find an element $\ob$ 
such that $(\cb,\ob)\in\cMC(\cH,\l,\n)$, 
$n'_{0,0}$ also vanishes and 
one gets a deformed (strict) $A_\infty$-algebra. 
Thus we obtain Theorem \ref{thm:defAinfty} above.

 \subsection{Gauge equivalence and moduli spaces}
\label{ssec:moduli}

Continuing with the general philosophy of deformation theory, 
we regard two deformations as equivalent if they are related 
by {\em gauge equivalence}, that is, 
if they differ by the action of the group obtained 
by exponentiating the action of $\g^0$ 
of the controlling dg Lie algebra $\g$. 

For the case of $L_\infty$-algebras instead of dg Lie algebras, 
it is more subtle to show that the gauge equivalence given in a 
similar way in fact defines an equivalence relation, 
that is, the composition of gauge transformations is a gauge
transformation. 
In order to avoid such conceptually irrelevant subtlety, 
we give a definition of gauge equivalence 
in a more formal way in terms of piecewise smooth paths, 
though these definitions should be equivalent under some appropriate 
assumptions (see \cite{FOOO}). 
\begin{defn}[Gauge equivalence]
Given an $L_\infty$-algebra $(\cH_c,\l)$, 
two elements $\cb_0\in\cMC(\cH_c,\l)$ and $\cb_1\in\cMC(\cH_c,\l)$ 
are called {\em gauge equivalent} 
iff there exists a piecewise smooth path 
$\cb_t\in\cMC(\cH_c,\l)$, $t\in [0,1]$ such that 
\begin{equation}
 \fd{t}\cb_t
=\sum_{k\ge 0}
\ov{k !}
 l_{1+k}(\alpha(t),\cb_t^{\otimes k})
 \label{diffeq-c}
\end{equation}
for a degree minus one element $\alpha(t)\in\cH_c^{-1}$. 
 \label{defn:gauge}
\end{defn}
By this definition, it is clear that the gauge equivalence actually 
defines an equivalence relation. 

One can also express this gauge transformation 
in terms of a path ordered integral 
as $c_1=c_1(\{l_k\},c_0,\alpha(t))=c_0+\cdots$ 
\cite{thesis,KaTe}. 
\begin{defn}[Moduli space]
For an $L_\infty$-algebra $(\cH_c,\l)$ 
and the solution space of its Maurer-Cartan equation $\cMC(\cH_c,\l)$, 
the corresponding moduli space $\cM(\cH_c,\l)$ is defined as 
\begin{equation*}
 \cM(\cH_c,\l):=\cMC(\cH_c,\l)/\sim\ ,
\end{equation*}
where $\sim$ is the gauge equivalence in Definition \ref{defn:gauge}. 

The moduli space for an $A_\infty$-algebra $(\cH_o,\m)$ 
is also defined in a similar way and denoted by 
$\cM(\cH_o,\m):=\cMC(\cH_o,\m)/\sim$. 
 \label{defn:moduli}
\end{defn}
The following classical fact is known 
(for instance see \cite{Ko1, Fukaya2, FOOO, thesis, KaTe}; 
some of these include the case of $A_\infty$-algebras, for which 
a similar fact holds). 
\begin{thm}
For two $L_\infty$-algebras $(\cH_c,\l)$ and $(\cH_c',\l')$, 
suppose there exists an $L_\infty$-morphism 
$\f: (\cH_c,\l)\raw (\cH_c',\l')$. 
Then there exists a well-defined map
\begin{equation*}
 \f_{\sim} : \cM(\cH_c,\l)\raw\cM(\cH_c',\l')
\end{equation*}
and in particular $\f_{\sim}$ gives an isomorphism if 
$\f$ is an $L_\infty$-quasi-isomorphism. 
 \label{thm:known}
\end{thm}
Then, as a corollary of Theorem \ref{thm:defs} we have the following: 
\begin{cor}[$A_\infty$-structure parameterized by the moduli space 
of $L_\infty$-structures]
For an $L_\infty$-algebra $(\cH_c,\l)$ and an $A_\infty$-algebra 
$(\cH_o,\m)$, 
suppose there exists 
an OCHA $(\cH=\cH_c\oplus\cH_o,\l,\n)$ such 
that $(\cH_o,\{n_{0,k}\})=(\cH_o,\m)$. 
Also, let $(\cH'_c,\l')$ be an $L_\infty$-algebra obtained by 
the suspension of the dg Lie algebra $\Coder(T^c A)$ 
with $D=[\m,\ ]$ and Lie bracket $[\ ,\ ]$. 
The OCHA $(\cH,\l,\n)$ then gives 
a map from $\cM(\cH_c,\l)$ to $\cM(\cH_c',\l')$ and 
it is in particular an isomorphism 
if the $L_\infty$-morphism $(\cH_c,\l)\raw(\cH_c',\l')$ is an 
$L_\infty$-quasi-isomorphism. 
 \label{cor:Ainfty}
\end{cor}
In a similar way as in the $A_\infty$ or $L_\infty$ case, 
we can define the moduli space of the solution space 
of the Maurer-Cartan equations for an OCHA. 
\begin{defn}[Open-closed gauge equivalence]
Given an OCHA $(\cH,\l,\n)$, 
we call two elements $(\cb_0,\ob_0)\in\cMC(\cH,\l,\n)$ and 
$(\cb_1,\ob_1)\in\cMC(\cH,\l,\n)$ 
{\em gauge equivalent} iff 
there exists a piecewise smooth path 
$(\cb_t,\ob_t)\in\cMC(\cH,\l,\n)$, $t\in [0,1]$ such that 
$\cb_t$ satisfies differential equation (\ref{diffeq-c}) and 
$\ob_t$ satisfies 
\begin{equation*}
 \fd{t}\ob_t
 =\sum_{p,q\ge 0}\ov{p !}
 n_{1+p,q}(\alpha(t),\cb_t^{\otimes p};\ob_t^{\otimes q}) 
+\sum_{p,q,q'\ge 0}\ov{p!}
 n_{p,q+1+q'}(\cb_t^{\otimes p};
\ob_t^{\otimes q},\beta(t),\ob_t^{\otimes q'}) 
\end{equation*}
for degree minus one elements 
$(\alpha(t),\beta(t))\in(\cH_c^{-1},\cH_o^{-1})$. 
 \label{defn:gaugeOC}
\end{defn}
By definition, when $(\cb_0,\ob_0)$ and 
$(\cb_1,\ob_1)$ are gauge equivalent in the sense of an OCHA, 
$\cb_0$ and $\cb_1$ are gauge equivalent in the
sense of the $L_\infty$-algebra. 
\begin{defn}[Moduli space for an OCHA]
For an OCHA $(\cH,\l,\n)$ 
and the solution space of its Maurer-Cartan equations $\cMC(\cH,\l,\n)$, 
the moduli space for the OCHA $(\cH,\l,\n)$ is defined by 
\begin{equation*}
 \cM(\cH,\l,\n):=\cMC(\cH,\l,\n)/\sim\ ,
\end{equation*}
where $\sim$ is the gauge equivalence in Definition \ref{defn:gaugeOC}.
\label{defn:moduli-oc}
\end{defn}
Then, due to the theorems in section \ref{sec:MMth} and 
in particular Corollary \ref{cor:inverse}, 
the following theorem is obtained in a similar way 
as in the $A_\infty$ and $L_\infty$-cases. 
\begin{thm}
Suppose we have an OCHA homomorphism 
$\f: (\cH,\l,\n)\raw (\cH',\l',\n')$ between two OCHAs. 
Then, $\f$ induces a well-defined map between two moduli spaces 
$\f_{\sim}:\cM(\cH,\l,\n)\raw \cM(\cH',\l',\n')$. 
Furthermore, if $\f$ is an OCHA quasi-isomorphism, 
it induces an isomorphism between the two moduli spaces. 
 \label{thm:moduli}
\end{thm}
Thus, the moduli space $\cM(\cH,\l,\n)$ is also a 
homotopy invariant notion and in particular 
the equivalence class of deformations given by 
Theorem \ref{thm:defAinfty} is described by $\cM(\cH,\l,\n)$.

\begin{center}
\noindent{\large \textbf{Acknowledgments}}
\end{center}

H.~K would like to thank A.~Kato, T.~Kimura, H.~Ohta, A.~Voronov 
for valuable discussions. 
Also, he is very grateful to 
the Department of Mathematics of the University of Pennsylvania 
for hospitality, where this work was almost completed. 
J.~S.  would like to thank 
E.~Harrelson, T. Lada, M. Markl and B. Zwiebach 
for valuable discussions.

\input{insert228.tex}


\end{document}

%% file: insert228.tex
\def\calP{{\mathcal P}} \def\fracg{{\frak g}}
\def\Leib{{{\mathcal L}\it eib}}
\def\black{\bullet}
\def\white{\circ} \def\ext{\mbox{\large$\land$}}

\def\fA{\frak{A}}
\def\fC{\frak{C}}
\def\a{\frak{a}}
\def\b{\frak{b}}
\def\X{\frak{X}}
\def\Y{\frak{Y}}

\section*{Appendix (by M.~Markl): 
Operadic interpretation of $A_\infty$-algebras
over $L_\infty$-algebras}

This part of the paper assumes some knowledge of the language of
operads and related notions, see the
book~\cite{MSS}, namely Section~II.3.7 of this
book.  We explain here how $A_\infty$-algebras over
$L_\infty$-algebras can be interpreted using a `colored' version of
the standard theory of strong homotopy algebras in the form formulated
in~\cite[Propostion~II.3.88]{MSS}. Let us
briefly recall some necessary background material.

Assume that $\calP$ is a quadratic Koszul operad governing the
algebraic structure we have in mind (such as associative algebra, Lie
algebra, etc.) and let $\calP^!$ denote the quadratic dual of $\calP$.
Proposition~II.3.88 of~\cite{MSS} then says
that a {\em strongly homotopy $\calP$-algebra\/} on a graded vector
space $V$ is the same as a degree $+1$ differential on the cofree
nilpotent $\calP^!$-coalgebra 
$T_{\calP^!}(\downarrow \hskip -.2em V)$
on the desuspension of $V$.  Using a colored version of this
proposition, we show that $A_\infty$-algebras over $L_\infty$-algebras
are in fact {\em strongly homotopy Leibniz pairs\/}.

Let $\rho : \fracg \to {\rm Der}\, A$ be a Leibniz pair as in
Definition~\ref{defn:g-alg}. These Leibniz pairs are algebras over a
two-colored operad $\Leib$, with the white color denoting
inputs/output in $\fracg$ and the black color inputs/output in
$A$. The operad $\Leib$ is a quadratic
$\{\white,\black\}$-colored operad generated by one antisymmetric
binary operation $l$ of type $(\white,\white) \to \white$ for the Lie
multiplication in $\fracg$, one binary operation 
$m$ of type
$(\black,\black) \to \black$ for the associative multiplication in
$A$, and one binary operation $\rho$ of type $(\white,\black) \to
\black$ for the action of $\fracg$ on $A$. The relations defining
$\Leib$ as a quadratic colored operad can be easily read off
from eq.(\ref{Xab})  and eq.(\ref{XYa}).
 
We may safely leave as an exercise to verify that the quadratic
dual $\Leib^!$ of the operad $\Leib$ describes objects 
$(\frak{C},\frak{A})$ 
consisting of a commutative associative algebra $\frak{C}$, 
an associative algebra $\frak{A}$ 
and an action of $\frak{C}$ on $\frak{A}$ that satisfies
\begin{equation*}
\X(\a\b) = \X(\a)\b = \a\X(\b),\ 
 \mbox { for $\X \in \fC$ and $\a,\b \in \fA$,}
\end{equation*}
and
\begin{equation*}
(\X\Y)(\a) = \X(\Y(\a)) = \Y(\X(\a)),\ 
 \mbox { for $\X,\Y \in \fC$ and $\a \in \fA$.}
\end{equation*}
It is equally simple to prove that the cofree nilpotent
$\Leib^!$-algebra cogenerated by a colored space 
$V_\white \oplus V_\black$ equals 
$C(V_\white) \otimes T^c(V_\black)$, where
$C(V_\white)$ is the cofree nilpotent cocommutative coassociative
coalgebra cogenerated by $V_\white$ and $T^c(V_\black)$ is the cofree
nilpotent coassociative coalgebra cogenerated by
$X_\black$ (the tensor coalgebra). 

It immediately follows from these calculations that the obvious
colored version of the above
mentioned~\cite[Propostion~II.3.88]{MSS}
identifies $A_\infty$-algebras over $L_\infty$-algebras in the
coalgebra description as in subsection~\ref{ssec:coalg} 
(where $V_\white$ and $V_\black$ correspond to $\cH_c$ and $\cH_o$) 
with strongly homotopy Leibniz pairs.

The only nontrivial thing which we left aside   
was to prove that $\Leib$ is a {\em Koszul\/} quadratic colored
operad, in the sense of~\cite{laan:03}. 
This is, according to~\cite{markl:haha},
necessary for the homotopy invariance of these 
strongly homotopy Leibniz pairs, 
though the above constructions make sense even without the Koszulity. 
We believe that Koszulness of $\Leib$ would follow from a spectral sequence
argument similar to that used in the proof of~\cite[Theorem~4.5]{markl:dl}.